\documentclass[final,times,authoryear,1p]{elsarticle}

\usepackage{amsmath,amsfonts,amssymb}
\usepackage{amsthm}
\usepackage{graphicx}
\usepackage{multirow}
\usepackage{algpseudocode}
\usepackage{algorithm}
\usepackage{color}

\usepackage{fancyhdr}
\newtheorem{theorem}{Theorem}[section]

\newtheorem{remark}[theorem]{Remark}

\usepackage{amssymb}

\journal{Transportation Research Part C}

\begin{document}

\begin{frontmatter}



\begin{center}
\textcolor{blue}{ARTICLE LINK:  http://www.sciencedirect.com/science/article/pii/S0968090X14003635
\\  PLEASE CITE THIS ARTICLE AS\\ 
Piccoli, B., Han, K., Friesz, T.L., Yao, T., Tang, J., 2015. Second order models and traffic data from mobile sensors. Transportation Research Part C 52, 32-56.}
 \line(1,0){469}
 \end{center}

\title{Second-order models and traffic data from mobile sensors}


\author[math1]{Benedetto Piccoli}
\ead{piccoli@camden.rutgers.edu}

\author[cts]{Ke Han\corref{cor}}
\ead{k.han@imperial.ac.uk}

\author[ie]{Terry L. Friesz}
\ead{tfriesz@psu.edu}

\author[ie]{Tao Yao}
\ead{tyy1@engr.psu.edu}

\author[cts]{Junqing Tang}
\ead{junqing.tang13@imperial.ac.uk}

\cortext[cor]{Corresponding author}

\address[math1]{Department of Mathematics, Rutgers University, Camden, NJ 08102, USA}
\address[cts]{Department of Civil and Environmental Engineering, Imperial College London, London SW7 2BU, UK}
\address[ie]{Department of Industrial and Manufacturing Engineering, Pennsylvania State University, University Park, PA 16802, USA}

\begin{abstract}
Mobile sensing enabled by GPS or smart phones has become an increasingly important source of traffic data. For sufficient coverage of the traffic stream, it is important to maintain a reasonable penetration rate of probe vehicles. From the standpoint of capturing higher-order traffic quantities such as acceleration/deceleration, emission and fuel consumption rates, it is  desirable to examine the impact on the estimation accuracy of sampling frequency on vehicle position. Of the two issues raised above, the latter is rarely studied in the literature.  This paper addresses the impact of both sampling frequency and penetration rate on mobile sensing of highway traffic. To capture inhomogeneous driving conditions and deviation of traffic from the equilibrium state, we employ the second-order {\it phase transition model} (PTM). Several data fusion schemes that incorporate vehicle trajectory data into the PTM are proposed. And, a case study of the NGSIM dataset is presented which shows the estimation results of various Eulerian and Lagrangian traffic quantities. The findings show that  while first-order traffic quantities can be accurately estimated even with a  low sampling frequency, higher-order traffic quantities, such as acceleration, deviation, and emission rate, tend to be misinterpreted due to insufficiently sampled vehicle locations. We also show that a correction factor approach has the potential to reduce the sensing error arising from low sampling frequency and penetration rate, making the estimation of higher-order quantities more robust against insufficient data coverage of the highway traffic.
\end{abstract}


\end{frontmatter}

\section{Introduction}

With the increased availability of mobile traffic data and the advancement of sensing technology, data collected through GPS, smart phones or other mobile devices have become a major source of traffic information  for various applications. Advantages of mobile sensing, in comparison with fixed-location sensing (e.g. using loop detectors and cameras), include potentially complete spatial and temporal coverage of traffic network and high positioning accuracy \citep{MC}.

Traffic data related to speed, density, queue size and travel time, which are categorized as {\it lower-order} quantities, can often be estimated in conjunction with first-order traffic flow models such as the Lighthill-Whitham-Richards (LWR) model \citep{LW, Richards} and the cell transmission model (CTM) \citep{CTM1}. \cite{SB} employ a weak formulation of boundary conditions for the LWR model based on the Godunov scheme, which is then applied to the I-80 highway dataset.  In \cite{CFP}, the LWR model is discretized in connection with initial and boundary conditions, which is applied to traffic estimation on a circular urban motorway using mobile data.   \cite{CC1} propose convex formulations for data assimilation using both Eulerian (fixed) and Lagrangian (mobile) traffic data based on the Hamilton-Jacobi representation of highway traffic and the generalized Lax-Hopf formula. \cite{WBTPB} employ a velocity-based LWR model with transformed fundamental diagram to perform data fusion using ensemble Kalman filter. Independently, \cite{YLHVS} consider the LWR model with a transformed coordinate system, namely the Lagrangian coordinates, and perform traffic estimation using extended Kalman filter. These studies mainly focus on freeway traffic.

In another line of research, mobile data have been used extensively for estimating queue size and delay at signalized intersections in arterial networks. \cite{BHHB} use sampled vehicle travel times  to estimate delay patterns near a signalized junction, where the authors use the LWR theory with a triangular fundamental diagram to express the relationship among flow, shock speed, queue size and queuing time. Following this work, \cite{BHS} devise a reverse modeling process that construct the dynamic queue length in real time. \cite{CQJR} further explore probe vehicle trajectories in estimating queue size in real time with the benefit of less communication cost in data collection. \cite{CC} propose a statistical approach for estimating queue length using an analytical formulation based on conditional probabilities. The authors also address the estimation accuracy with a wide range of probe penetration rates. \cite{ACXS} consider several {\it measures of effectiveness} and estimation methods to identify proper penetration rates.

Although the LWR model and the CTM have been used effectively in estimating lower-order quantities, they have been used less frequently in estimating high-order traffic quantities such as acceleration/deceleration, deviation (perturbation), emission and fuel consumption rates. There exist a number of attempts to estimate acceleration/deceleration or emission rates through differentiating macroscopic traffic quantities analytically or numerically (e.g. \cite{LKVZSV}). However, in this process higher-order variations inherent in these quantities, typically on a microscopic scale, are insufficiently captured due to the low temporal-spatial resolution of the traffic data and the discrete models. To fill this gap, this paper proposes a second-order traffic flow model supported by high-resolution mobile data to address the issue of estimating high-order traffic quantities. Unlike most existing studies on mobile sensing which primarily focus on probe penetration rate \citep{DLWLW, KMPV, YC}, we consider the additional effect of under sampling on the estimation accuracy. This is a concern because most mobile data provide location or speed of a moving vehicle every 3-4 seconds (such as GPS data), but higher-order variations in speed, acceleration and emission may take place on a smaller time scale; this is true especially for congested and unstable traffic.  Such an observation has raised the need to examine the efficacy of existing sensing technologies and estimation methods in reconstructing the profiles for these higher-order quantities.

The second-order traffic flow model employed by this paper is the hyperbolic {\it phase transition model} (PTM),  first introduced by \cite{Colomboa} and studied subsequently by \cite{Colombo and Corli, CGP} and \cite{BABW}. Second- and higher-order traffic models were proposed by many researchers to overcome some limitations of the LWR model in describing complex waves observed in vehicular traffic. Most second-order models tend to pose, in addition to the LWR-type equation for the conservation of vehicles, a second equation for the conservation or balance of momentum. One of the first such models is by Payne and Whitham in the 1970s
\citep{Payne1,Whitham}. In a celebrated paper \citep{Daganzo1995} Daganzo criticized second-order models by 
showing various drawbacks including the possibility of cars going backward.
Most of such drawbacks were later addressed by the Aw-Rascle-Zhang model,
independently proposed by  \cite{AR}
and \cite{Zhang}.
More recently, the phase transition models drew increased attention from researchers
for their capability of representing complex waves while keeping
the LWR structure for light traffic, see  \cite{BWGPB} for more details.

Let us further comment on the possible use of second-order model for urban arterial traffic. The aforementioned complex wave phenomenon, well captured by
second-order models, are mainly observed in highway traffic.
Indeed, phantom waves, stop-and-go waves and others need a
long stretch of road with no interruption to manifest themselves.
The situation of arterial traffic is quite different because of the
presence of many junctions with traffic signs or signals.
Since the LWR model captures well backward wave propagation from junctions
or signals, and vehicle movements on arterials are relatively uniform, it is typically sufficient to describe arterial traffic, while very limited to model more complex waves in highway traffic.

In order to address issues related to higher-order traffic quantities, we focus on the  {\it Next Generation SIMulation} (NGSIM) dataset \footnote{http://ngsim-community.org/}. The NGSIM program collected high-quality traffic and vehicle trajectory data on a stretch of I-80 highway in California. A total of 45 minutes of traffic data were collected, segmented into three 15 minute periods.  The dataset contains vehicle trajectory recorded at a high resolution of every 0.1 s. Derived information on instantaneous velocity and acceleration is also available. A detailed description of the NGSIM field experiment and data will be provided in Section \ref{secNumerical}.

\subsection{Contribution and findings}

This paper addresses the effectiveness of mobile sensing and traffic estimation from two aspects: sampling frequency and probe penetration rate. The subject of estimation will include both lower-  and higher-order traffic quantities. To support this study, we employ the second-order phase transition model (PTM) as well as its modifications. The underpinning dataset is provided by the NGSIM program, which contains high resolution vehicle trajectories on a segment of I-80.  Such detailed vehicle trajectories provide unique information regarding the higher-order variations of traffic, which is unavailable through traditional mobile or fixed sensors.  We propose several data fusion schemes that integrate mobile data into the PTM under various assumptions. As a result, vehicle speed, acceleration, density, deviation (perturbation),  engine power demand and emission rate  can be estimated along vehicle trajectories (these quantities are Lagrangian as they are associated with moving vehicles). We demonstrate how these estimations deteriorate with less frequent sampling of vehicle locations, with special attention given to the difference between lower- and higher-order quantities.  We also present a method for estimating Eulerian traffic quantities; i.e. quantities associated with a give point in the temporal-spatial domain. The joint effect of sampling frequency and probe penetration rate in reconstructing lower- and higher-order quantities are assessed. It is found that, for both Lagrangian and Eulerian estimations, lower-order traffic variables (such as density) can be estimated relatively accurately, while higher-order ones (such as emission and perturbation) are far more susceptible to under sampling and lower penetration rate. This is consistent with the expectation that higher-order variations tend to be overlooked by data presented at a low resolution. We also show that Eulerian estimation of higher-order quantities enjoy improved accuracy than Lagrangian estimation due to the effect of spatial-temporal aggregation. And, the efficacy of Eulerian estimation can be further improved by applying a correction factor approach. This is demonstrated by using hydrocarbon emission rate as an example and proven to be an effective estimation method even with low sampling frequency and penetration rate.

The main contributions and  findings made by this paper are summarized as follows. 

\begin{itemize}
\item Several methods for integrating vehicle trajectory data into the phase transition model are proposed. These methods are further extended to perform both Lagrangian and Eulerian estimation of traffic quantities.

\item Along each vehicle trajectory, the estimation of higher-order traffic quantities deteriorates significantly when the sampling frequency decreases, while the estimation of first-order quantities remain relatively accurate with the same sampling frequency.

\item For the Eulerian estimation, we provide numerical results on density (first-order) and HC emission rate (higher-order) with varying sampling frequency and penetration rate.  The HC estimation enjoys an improved accuracy compared to the Lagrangian estimation due to the averaging of multiple measurements, although it also deteriorates with under sampling and lower penetration rate.

\item We provide a comparative study of first-order (LWR) models and the phase transition model in terms of estimating first-order density (notice that higher-order quantities such as acceleration and perturbation are not captured explicitly by the LWR model). The PTM outperforms the LWR model in terms of estimation accuracy. In addition,  deviation of individual vehicles from the equilibrium state, which is not captured by the LWR model, has been illustrated.

\item A correction factor approach is proposed for estimating HC emission rate with insufficient mobile data coverage. We employ a regression model to correct the predicted emission rate in order to minimize the misinterpretation of the ground truth caused by lower sampling frequencies or probe penetration rates. The validity and effectiveness of this approach is shown through cross validation. 
\end{itemize}

The rest of this paper is organized as follows. Section \ref{secMeasurement} recaps several numerical approximations of basic traffic quantities based on vehicle trajectory data. Section \ref{secTrafficModel} introduces the phase transition model (PTM), followed by three data fusion schemes. Models for HC emission rate and engine power demand are also described therein. Section \ref{secsmooth} illustrates the procedure employed to reduce random noises in the raw dataset and errors derived from numerical differentiation.  In Section \ref{secNumerical}, we calibrate the PTM using the NGSIM data. Section \ref{secestimatetraffic} assesses the estimation quality for first- and second-order quantities along vehicle trajectories, when the sampling frequency varies. Section \ref{secEulerian} and \ref{secHC} perform Eulerian estimation for density and HC emission rate respectively, and evaluate its effectiveness against various sampling frequencies and penetration rates. A correction factor approach approach for improving the accuracy and reliability of HC emission is also presented and tested. Section \ref{secconclusion} provides some concluding remarks.

\section{Estimating basic traffic quantities}\label{secMeasurement}

Onboard sensors such as GPS measure the location of a moving car every $\delta t$ seconds, where $\delta t$ is related to the device's characteristics such as desired precision and transmission capacity. For a given vehicle, we denote by $x(t)$ its location and by $v(t)$ its velocity at time $t$. Assume that the location is recorded at three consecutive time stamps $t_1, t_2$ and $t_3$ with $t_2-t_1=t_3-t_2=\delta t$. From these measurements one can deduce the approximate velocities in the time intervals $[t_1,t_2]$ and $[t_2,t_3]$ respectively as
\begin{equation}\label{approxv}
v_{1,2}~=~\frac{x(t_2)-x(t_1)}{\delta t},\qquad\qquad
v_{2,3}~=~\frac{x(t_3)-x(t_2)}{\delta t}
\end{equation}
The velocity at time $t_2$ is approximated as
\begin{equation}\label{approxmv}
v(t_2)~\approx~{v_{1,2}+v_{2,3}\over 2}~=~{x(t_3)-x(t_1)\over 2\delta t}
\end{equation}
One also gets estimate for the acceleration:
\begin{equation}\label{approxa}
a(t_2)~=~\frac{D}{Dt}v(t,\,x)\Big\vert_{t=t_2}~\approx~ \frac{v_{2,3}-v_{1,2}}{\delta t}~=~{x(t_3)-2x(t_2)+x(t_1)\over \delta t^2}
\end{equation}
where $v(t,\,x)$ denotes the velocity function in the Eulerian coordinates; $D/Dt=d/dt+v\cdot d/dx$ represents the material derivative in the Eulerian coordinates corresponding to the acceleration of the car in the Lagrangian coordinate.

Another important quantity to estimate is the spatial variation of velocity  in the Eulerian coordinate.  For notation convenience, we set $x_i\doteq x(t_i)$. Assuming a mild variation in time of the Eulerian velocity $v(t,x)$,  we write:
\[
v\left(t_2,\,\frac{x_2+x_1}{2}\right)\approx v_{1,2},\qquad\qquad
v\left(t_2,\,\frac{x_3+x_2}{2}\right)\approx v_{2,3}
\]
from which we get, by setting $\delta x= \frac{x_3+x_2}{2}-\frac{x_2+x_1}{2}$, that
\[
\frac{\partial}{\partial x}v\left(t_2,~\frac{\frac{x_3+x_2}{2}+\frac{x_2+x_1}{2}}{2}\right)\approx \frac{v_{2,3}-v_{1,2}}{\delta x}
\]
In other words,
\begin{equation}\label{approxvx}
\frac{\partial}{\partial x}v\left(t_2,\frac{x_3+2x_2+x_1}{4}\right)\approx \frac{v_{2,3}-v_{1,2}}{\frac{x_3-x_1}{2}}~=~{2\over \delta t}{x_3-2x_2+x_1\over x_3-x_1}
\end{equation}
Clearly such approximation is acceptable as long as the variation between $v_{1,2}$
and $v_{2,3}$ is not too large.

\section{Model fitting using vehicle trajectory data}\label{secTrafficModel}
We consider the phase transition model (PTM) as well as some of its variations to represent the dynamics of highway traffic.  In addition, we also consider models proposed by \cite{Modal} and \cite{Ahn} to describe instantaneous engine power demand and emission rate based on detailed vehicle trajectory. All these models will be detailed in this section.

\subsection{The phase transition model}\label{secPTM}

The hyperbolic {\it phase transition model} (PTM) for traffic flow belongs to a class known as {\it second-order}, since it captures second-order variations of traffic in addition to the average velocity and density. Other second-order models include the Payne-Whitham model proposed independently by \cite{Payne1,Payne2} and \cite{Whitham}, and the Aw-Rascle-Zhang model developed by \cite{AR} and \cite{Zhang}. The phase transition model is motivated by the empirical observation that when the vehicle density exceeds certain critical value, the density-flow pairs  are scattered  in a two-dimensional region, instead of forming a one-to-one relationship. This is in contrast to the first-order kinematic wave models such as the classical Lighthill-Whitham-Richards (LWR) model \citep{LW, Richards}.

The phase transition model considers of two distinct traffic phases: the uncongested phase and the congested phase. In the uncongested phase, the dynamic is governed by the LWR model
\begin{equation}\label{PTMfree}
\begin{cases}
\partial_t   \rho(t,x)+\partial_x \left[\rho(t,x)\cdot v(t,x)\right]~=~0\\
v(t,x)~=~v\big(\rho(t,x)\big)
\end{cases}
\end{equation}
where $\rho(t,\,x)$ denotes density, and the velocity $v\big(\rho(t,\,x)\big)$ is expressed as  a function of density. The product $\rho(t,x)\cdot v(t,x)$ represents flow. On the other hand, the congested phase is governed by the following system of conservation laws:
\begin{equation}\label{PTMcongested}
\begin{cases}
\partial_t \rho(t,x)+\partial_x\left[\rho(t,x)\cdot v(t,x)\right]~=~0\\
\partial_t  q(t,x)+ \partial_x \left[(q(t,x)-q^*)\cdot v(t,x)\right]~=~0\\
v(t,x)~=~v\big(\rho(t,x),\,q(t,x)\big)
\end{cases}
\end{equation}
where the velocity $v(\rho,\,q)$ depends not only on the local density $\rho$, but also on $q$ which describes the perturbation or deviation from the equilibrium state. $q^*$ is a given constant corresponding to the equilibrium (stationary) traffic state. Some choices of the functional form $v(\rho,\,q)$ include
\begin{equation}\label{funcv1}
v(\rho,\,q)~=~\left(1-{\rho\over\rho_{jam}}\right)\cdot {q\over \rho}\qquad\hbox{or}\qquad v(\rho,\,q)~=~v(\rho)(1+q)
\end{equation}
where $\rho_{jam}$ denotes the jam density, $v(\rho)$ denotes the equilibrium velocity function. In this paper, we employ the following functional form for the velocity.
\begin{equation}\label{funcv2}
v(\rho,\,q)~=~A(\rho_{jam}-\rho)+B(q-q^*)(\rho_{jam}-\rho)
\end{equation}
for some positive parameters $A$, $B$ and $\rho_{jam}$, where $q-q^*\in[-1,\,1]$ \footnote{In general, the deviation $q-q^*\in[q_{min},\,q_{max}]$ for some lower and upper bounds $q_{min}<0$ and $q_{max}>0$. This paper assumes that the deviations are symmetric around the equilibrium state; that is $q-q^*\in[-q_{max},\,q_{max}]$ for some $q_{max}>0$. Without loss of generality, we further assume that $q-q^*\in[-1,\,1]$ since $q_{max}$ can be factored into the coefficient $B$, which will later be calibrated using field data in Section \ref{secNumerical}.}.  
\begin{remark}
Several factors should be taken into account when formulating an appropriate functional form for $v(\rho,\,q)$: 
\begin{itemize}
\item[(i)] When $q=q^*$ (the equilibrium state), the velocity $v(\rho,\,q)$ coincides with the stationary density-velocity relationship.

\item[(ii)] $v(\rho,\,q)\geq 0$, and zero velocity arises only at densities corresponding to the jam density in the equilibrium state (this property is posed by \cite{GP2009} to ensure physically correct behavior of models on networks).

\item[(iii)] Assume that $q-q^*\in[q_{min},\,q_{max}]$. The 2-D region formed by the curves $v(\rho,\,q^*+q_{min})$  and $v(\rho,\,q^*+q_{max})$  should match the empirical  density-velocity data plot (see Figure \ref{figfdsall} for an example).
\end{itemize} 
As we show later in the model calibration (Section \ref{secNumerical}), the density-velocity plot suggests a well-defined affine stationary relation for the congested branch, which is represented by the first term of \eqref{funcv2}. In addition, the upper and lower envelops of the 2-D region, expressed respectively as $(A+B)(\rho_{jam}-\rho)$ and $(A-B)(\rho_{jam}-\rho)$ also match well with the data points; see Figure \ref{figfdsall}. Finally, the proposed function \eqref{funcv2} clearly satisfies requirements (i) and (ii). This justifies our  choice of $v(\rho,\,q)$.
\end{remark}

For the PTM in general, there is some flexibility in choosing the system of equations for the congested phase. For instance, \cite{Gotin} proposes a phase transition model which employs the Aw-Rascle-Zhang equations \citep{AR, Zhang} for the congested phase. Furthermore, one may consider a version of the PTM by taking into account the reaction time of drivers. More specifically, following \cite{SM}, we write 
\begin{equation}\label{PTMsm}
\partial_t q(t,x)+\partial_x[\big(q(t,x)-q^*\big)\cdot v(t,x)]~=~{q^*-q(t,x)\over T-\tau}
\end{equation}
as the second equation for the congested phase. The non-vanishing right hand side is called the Siebel-Mauser source term. Here, $\tau$ (in second) is a reaction time and typically varies within $[0.5,\,2]$ \citep{Koppa}. In this paper, the value of $\tau$ is set to be $1$ (second). Following \cite{SM} and \cite{KM}, we choose $T=2/3$ (second) and therefore $T-\tau=-1/3$. Finally $T-\tau$ is modeled as a factor which can be  positive (for very small or very high densities) or negative (for intermediate densities). One should realize  that a negative factor gives rise to stable traffic, while a positive ones produces instabilities \citep{SM}.

\subsection{Estimating traffic quantities using the PTM}

The primary variables of the PTM are vehicle density $\rho$ (first-order) and deviation  $q$ (second-order). It is possible to estimate both quantities within the PTM using just vehicle trajectory data, as we demonstrate below. Note that subsequent derivations are established based on the congested phase of the PTM for two reasons: (1) the NGSIM dataset suggests a prevailing congested condition of traffic during the time of study; (2) it is mainly within the congested phase where traffic instability and heterogeneity driving conditions arise. 

We propose three methods to perform the estimation of the two variables mentioned above. These methods rely on different assumptions made on the model and the traffic conditions.

\subsubsection{Phase transition model with source term: Method 1}
We consider the  PTM  with velocity in the congested phase expressed as
\begin{equation}\label{ptmv}
v(\rho,q)=A(\rho_{jam}-\rho)+B(q-q^*)(\rho_{jam}-\rho)
\end{equation}
In addition, we add a Siebel-Mauser type
source term. Thus the equation for the congested phase becomes
\begin{equation}\label{eqn1}
q_t+[(q-q^\ast)\cdot v]_x=\frac{q^*-q}{T-\tau}
\end{equation}
where, without causing any confusion, we simplify the notations by dropping $(t, x)$ and using subscripts $t$ and $x$ to denote partial derivatives with respect to time and location. Such notation convention will be adopted throughout the rest of this paper.

Method 1 employs the following simplifications. \\

\noindent {\bf (Assumption 1).} \textit{Assume small variations of $\rho$ in $t$, and $\rho$, $q$ in $x$. That is, $\rho_t,\rho_x,q_x\approx 0$ with only $q_t$ being non-negligible.}  \\

\noindent By {\bf Assumption 1}, \eqref{eqn1} reduces to 
\begin{equation}\label{strongstability}
q_t~=~{q^*-q\over T-\tau}
\end{equation}
Then we can derive
\begin{align}
\frac{D}{Dt}v&~=~\partial_tv\big(\rho(t,\,x),\,q(t,\,x)\big)+v\cdot\partial_xv\big(\rho(t,\,x),\,q(t,\,x)\big) \nonumber\\
&~=~v_{\rho}\cdot\rho_t+v_q\cdot q_t+v\cdot \big(v_{\rho}\cdot \rho_x+v_q\cdot q_x\big)\nonumber\\
&~\approx~v_q\cdot q_t\nonumber\\
&~=~B(\rho_{jam}-\rho)\cdot {q^*-q\over T-\tau}\nonumber\\
\label{eqn2}
&~=~{1\over T-\tau}B(\rho_{jam}-\rho)(q^*-q)~=~{1\over T-\tau} \big(A(\rho_{jam}-\rho)-v\big) 
\end{align}
Taking into account only the measurements of $v$ and $Dv/Dt$, and by introducing variables $\hat\rho\doteq \rho_{jam}-\rho$, $\hat q\doteq q-q^*$, we deduce from \eqref{eqn2} that
\begin{align}\label{eqn3}
\hat\rho~=~&(\rho_{jam}-\rho)~=~\frac{1}{A}\left(v+(T-\tau)\frac{Dv}{Dt}\right)
\\
\label{eqn4}
\hat q~=~&(q-q^*)~=~-{A(T-\tau)\over B}\frac{\frac{Dv}{Dt}}{v+(T-\tau)\frac{Dv}{Dt}}
\end{align}
Following our discussion at the end of Section \ref{secPTM}, we take $T-\tau=-{1\over 3}$. The velocity $v$ is estimated according to \eqref{approxmv}, and the acceleration $Dv/Dt$ is computed from \eqref{approxa}.

\subsubsection{Phase transition model with source term: Method 2}

The second method again assumes a Siebel-Mauser source term, and relies on the following assumption. \\

\noindent {\bf (Assumption 2).} \textit{Assume that  $\rho_t$, $q_x\approx 0$}. \\

\noindent Then equation \eqref{eqn1} becomes
\begin{equation}\label{eqn5}
q_t+v_x(q-q^*)~=~{q^*-q\over T-\tau}
\end{equation}
We can now write
\begin{equation}\label{eqn6}
\frac{Dv}{Dt}~=~v_t+vv_x\approx v_qq_t+vv_x~=~v_q  \left({q^*-q\over T-\tau}-(q-q^*)v_x\right)+v v_x
\end{equation}
Recalling the variables $\hat\rho \doteq \rho_{jam}-\rho,~\hat q\doteq q-q^*$, we  obtain
\begin{equation}\label{eqn7}
v~=~A(\rho_{jam}-\rho)+B(q-q^*)(\rho_{jam}-\rho)~=~\hat\rho\,(A+B\hat q)
\end{equation}
and deduce from \eqref{eqn6}-\eqref{eqn7} that
\begin{equation}\label{eqn8}
\frac{Dv}{Dt}-vv_x~=~-v_q\left({\hat q\over T-\tau} + \hat q v_x\right)~=~-B\hat\rho\left(\frac{\hat q}{T-\tau}+\hat q v_x\right)
\end{equation}
From \eqref{eqn7} and  \eqref{eqn8} we immediately get the expressions for $\hat\rho$ and $\hat q$ in terms of $v,\,v_x$ and ${Dv\over Dt}$.
\begin{align}\label{eqn9}
\hat \rho~=~&{1\over A}{v + (T-\tau){Dv\over Dt}\over 1+ (T-\tau)v_x}
\\
\label{eqn10}
\hat q~=~&{A\over B}{(T-\tau)\left(vv_x-{Dv\over Dt}\right)\over v + (T-\tau){Dv\over Dt}}
\end{align}

\noindent The quantities ${Dv\over Dt},\,v_x$ are given by \eqref{approxa} and \eqref{approxvx} respectively.  Regarding velocity $v$, notice that if one approximates $v$ with \eqref{approxv}, then 
\begin{equation}\label{vanishing}
{Dv\over Dt}-vv_x~\approx~{x(t_{i+1})-2x(t_i)+x(t_{i-1})\over \delta t^2}-{x(t_{i+1})-x(t_{i-1})\over 2\delta t}\cdot {2\over \delta t}\cdot {x(t_{i+1})-2x(t_i)+x(t_{i-1})\over x(t_{i+1})-x(t_{i-1})}~\equiv~0
\end{equation}
which renders \eqref{eqn10} identically zero.  To avoid such a trivial case, one should instead approximate $v$ by
\begin{equation}\label{modifieddiff}
v(t_2)~\approx~ v_{1,2}~=~{x(t_2)-x(t_1)\over \delta t},\qquad \hbox{o}r\qquad  v(t_2)~\approx~ v_{2,3}~=~{x(t_3)-x(t_2)\over \delta t},
\end{equation}
or use robust numerical differentiation illustrated later in Section \ref{secrobustdiff}.

\subsubsection{Phase transition model without source term: Method 3}\label{secPTMnosource}
We consider the original PTM without the Siebel-Mauser type source term:
\begin{equation}\label{eqn11}
\begin{cases}
\rho_t+(\rho\cdot v)_x~=~0\\
q_t+((q-q^*)\cdot v)_x~=~0\\
v(\rho,\,q)~=~A(\rho_{jam}-\rho)+B(q-q^*)(\rho_{jam}-\rho)
\end{cases}
\end{equation}
Computation in this case is less straightforward than the previous two as no simplification is made in this case. We start with the identity
\[
\frac{Dv}{Dt}=v_t+v\cdot v_x=v_\rho\,\rho_t+v_q\,q_t+v(v_\rho\,\rho_x+v_q\,q_x)
\]
Using \eqref{eqn11}, we have
\begin{align}
{Dv\over Dt}&~=~v_\rho(-\rho_xv-\rho v_x)+v_q(-q_xv-(q-q^*)v_x)+v(v_\rho\rho_x+v_qq_x)\nonumber  \\
&~=~-v_x\big(v_\rho\rho+v_q(q-q^\ast)\big)\nonumber\\
\label{eqn12}
&~=~v_x\big(A\rho+B(q-q^\ast)(2\rho-\rho_{jam})\big)
\end{align}

\noindent Recall the variables $\hat\rho=\rho_{jam}-\rho$ and $\hat q=q-q^*$. 
Combining \eqref{eqn12} with the expression of $v(\rho,\,q)$ and solving for $\rho$ we get
\begin{equation}\label{eqn14}
\rho=\frac{1}{B\hat q}\left(v+{1\over v_x}\frac{Dv}{Dt}-A\rho_{jam}\right)
\end{equation}
One immediate observation from \eqref{eqn14} is that $\hat q$ and $v+{1\over v_x}{Dv\over Dt}-A\rho_{jam}$ always have the same sign. Substituting \eqref{eqn14} into \eqref{eqn12}, we get
\begin{equation}\label{quadratic}
B^2\rho_{jam}\hat q^2+B\left(2A\rho_{jam}-{1\over v_x}{Dv\over Dt}-2v\right)\hat q -A\left(v+{1\over v_x}{Dv\over Dt}-A\rho_{jam}\right)~=~0
\end{equation}
which is a quadratic equation in the variable $\hat q$. The discriminant of such quadratic equation is 
\begin{align}
\Delta~=~&4B^2\left(v+{1\over 2}{1\over v_x}{Dv\over Dt}-A\rho_{jam}\right)^2+4AB^2\rho_{jam}\left(v+{1\over v_x}{Dv\over Dt}-A\rho_{jam}\right)\nonumber \\
\label{discriminant}
~=~&4B^2\left(\left(v+{1\over 2}{1\over v_x}{Dv\over Dt}\right)^2-A\rho_{jam}v\right)
\end{align}
In order for any meaningful real root of \eqref{quadratic} to exist, a necessary condition is that $\Delta$ is nonnegative, that is,
\begin{equation}\label{necessary}
\left|v+{1\over 2}{1\over v_x}{Dv\over Dt}\right|~\geq~\sqrt{A\rho_{jam}v}
\end{equation}
If the equality holds in \eqref{necessary}, a real solution of $\hat q$ exists if and only if 
\begin{equation}
v+{1\over v_x}{Dv\over Dt}-A\rho_{jam}~<~0
\end{equation}

\noindent In the case where the strict inequality of \eqref{necessary} holds,  two distinct roots $\hat q_1$ and $\hat q_2$ exist. We distinguish between two cases:
\begin{itemize}
\item If 
\begin{equation}\label{nscase1}
v+{1\over v_x}{Dv\over Dt}-A\rho_{jam}~>~0
\end{equation} 
Then $\hat q_1\hat q_2<0$.  Equation  \eqref{quadratic} has one positive root and one negative root. By \eqref{eqn14}, one should choose the positive root since $\rho$ must be non-negative.
\item If 
\begin{equation}\label{nscase2}
v+{1\over v_x}{Dv\over Dt}-A\rho_{jam}~\leq~0
\end{equation}
Then $\hat q_1\hat q_2\geq 0$. Equation \eqref{quadratic} has two roots with the same sign\footnote{zero is considered to have both positive and negative signs.}. In view of \eqref{eqn14} and \eqref{nscase2}, to ensure that $\rho$ is nonnegative, both $\hat q_1$ and $\hat q_2$ should be nonpositive and at least one root is negative. This in turn requires that 
\begin{equation}\label{nseqn1}
2A\rho_{jam}-{1\over v_x}{Dv\over Dt}-2v~>~0
\end{equation}
In view of \eqref{nscase2}, a sufficient condition for \eqref{nseqn1} to hold is ${1\over v_x}{Dv\over Dt}>0$.
\end{itemize}

\noindent It turns out that the above analysis can be more easily interpreted in a discrete-time setting. Let us consider three consecutive time points $t_{i-1},\,t_{i},\,t_{i+1}$. First, notice that 
\begin{align}
{1\over v_x}{Dv\over Dt}~\approx~&{\delta t\over 2}\cdot{x(t_{i+1})-x(t_{i-1})\over x(t_{i+1})-2x(t_i)+x(t_{i-1})}\cdot {x(t_{i+1})-2x(t_i)+x(t_{i-1})\over \delta t^2} \nonumber
\\
\label{equlv}
~=~&{x(t_{i+1})-x(t_{i-1})\over 2\delta t}~\approx~v(t_i)~\geq~0
\end{align}
In light of this calculation, the feasibility condition \eqref{necessary} becomes
\begin{equation}\label{necessarydt}
{3\over 2}v~\geq~\sqrt{A\rho_{jam}v}\,,\qquad\hbox{or}\qquad v~\geq~{4\over 9}A\rho_{jam}
\end{equation}

\noindent Using  \eqref{equlv}, the decision rules following \eqref{necessary} can be explicitly summarized as follows.

\begin{algorithm}          
\caption{PTM without source term}
\label{alg1}            
\begin{algorithmic}                   
\If {$\displaystyle v~<~{4\over 9}A\rho_{jam}$\quad} \quad {the model has no solution}

\ElsIf  {$\displaystyle v~=~{4\over 9}A\rho_{jam}$\quad} \quad {$\displaystyle\rho={\rho_j\over 3} ,\quad \hat q=-{A\over 3 B}$} 

\Else 

         \If {$\displaystyle v~>~{1\over 2}A\rho_{jam}$\quad} \quad {$$ \rho= {2\rho_{jam}(2v-A\rho_{jam})\over 3v-2A\rho_{jam}+\sqrt{9v^2-4A\rho_{jam}v}}\,,\qquad\hat q~=~{3v-2A\rho_{jam}+\sqrt{9v^2-4A\rho_{jam}v}\over 2B\rho_{jam}}$$}
         
         \Else{} $$\rho={2\rho_{jam}(2v-A\rho_{jam})\over 3v-2A\rho_{jam}-\sqrt{9v^2-4A\rho_{jam}v}}\,, \qquad \hat q= {3v-2A\rho_{jam}-\sqrt{9v^2-4A\rho_{jam}v}\over 2B\rho_{jam}}$$
         
         \EndIf

\EndIf
\end{algorithmic}
\end{algorithm}

\begin{remark}
We note that the above computational procedure does not produce any result if the velocity is below ${4\over 9}A\rho_{jam}$.   Such an assumption may be rather restrictive in application, especially for traffic where low speed is prevailing, although this method is exact and consistent with the original phase transition model. As we later show in Table \ref{tabssrhoqentireperiod}, Method 3 is more resilient to under sampling in estimating $\hat q$ than the other two methods, whenever applicable.
\end{remark}

\subsection{Estimating power demand and vehicle emission rates}\label{secEmission}

Car emission and/or fuel consumption rates are considered in this section. They belong to high-order quantities as they are usually expressed as nonlinear functions of instantaneous speed and acceleration \citep{Post, Modal}. That is, the emission and/or fuel consumption rates are related to the modal operation of a moving vehicle such as idle, steady-state cruise, acceleration or deceleration. It requires detailed vehicle trajectory data in order to capture these driving conditions.

It has been shown that  fuel consumption is closely related to the engine power expressed as the {\it power demand function} (PDF) \citep{Post, Modal}, which in turn can be expressed as a nonlinear function of vehicle speed and acceleration, along with parameters related to road grade, air drag and other physical conditions. One example of such a PDF is provided in \cite{Modal} as
\begin{equation}\label{ztot}
Z~=~{M\over 1000}\cdot v\cdot(a+g\sin\theta)+\left(M g C_r+ {\sigma\over 2} v^2 A_c C_d\right)\cdot{v\over 1000}
\end{equation}
where $Z$ (in kw) is the power demand; $v$ (in km/h) and $a$ (in km/h/s) denote instantaneous speed and acceleration, respectively. $M$ (in kg) is vehicle mass and is chosen uniformly as 1200 in this paper for simplicity; $g=9.81 (m/s^2)$ is the gravitational constant; $\theta$ is the road grade and is chosen to be zero; $\sigma=1.225$ ($kg/m^3$) is the mass density of air.  Finally, $C_r=0.005$, $A_c=2.6$ ($m^2$), and $C_d=0.3$ denote respectively the rolling resistance coefficient, the cross-sectional area, and the aerodynamic drag coefficient, and their values are suggested by \cite{Chan}.

Regarding vehicle emission, we consider hydrocarbon (HC) and employ the model proposed in \cite{Ahn}. This paper presents a hybrid $3^{rd}$-order polynomial regression model with logarithmic transformation where the acceleration case ($(a\geq 0$) and the deceleration case ($a<0$) are fitted separately. The HC emission rate is modeled as 
\begin{equation}\label{hcrate}
\ln r^{HC}~=~
\begin{cases}
\displaystyle \sum_{i=0}^3\sum_{j=0}^3 U_{ij} \cdot v^i\cdot a^j \qquad\hbox{if }~ a~\geq~0
\\
\displaystyle \sum_{i=0}^3\sum_{j=0}^3 V_{ij} \cdot v^i\cdot a^j \qquad\hbox{if }~ a~<~0
\end{cases}
\end{equation}
where the matrices $\{U_{ij}\}_{i,j=0}^3$ and $\{V_{ij}\}_{i,j=0}^3$ are explicitly given in \cite{Ahn}, and are omitted here. This model is shown to better approximate the HC emission for combinations of relatively high speed and acceleration, and are thus chosen here since part of the monitored traffic is highly unstable containing accelerations up to 10 km/h per second.

\begin{remark}
Here we assume, as in almost all macroscopic traffic models, that all vehicles have the same characteristics, e.g. weight, engine type, size, and so on. The NGSIM dataset contains vehicle type which could potentially improve the estimation accuracy. This is, however, not pursued by this paper as we aim to describe the evolution of traffic at an aggregate level with bounded granularity on the characteristics of the flows. 
\end{remark}

We further note the fact that other types of models for instantaneous emission rate, power demand function and/or fuel consumption rate can be equally applied to our framework, as long as they are expressible as functions of speed, acceleration, and vehicle-specific characteristics. Although employing these alternative models will yield quantitatively different results, the general conclusion drawn from this paper will remain valid; that is, under sampling of the vehicle location and lower probe penetration rate tend to misinterpret these quantities (see Section \ref{secNumericalp} and \ref{subsecbinemission}), and the degree of such miscalculation depends on a number of factors, e.g. traffic condition (stable/unstable), vehicle type, and model parameters.

\section{Reconstructing ground-truth data}\label{secsmooth}

The NGSIM contains raw data  in  video format on vehicle trajectories with a high time resolution. Derived vehicle-specific quantities such as velocity and acceleration are also included therein. However, sensing errors are expected to exist, as are derived errors arising from numerical differentiation.  There exist a number of ways to pre-process those vehicle trajectory data via data smoothing techniques and robust numerical differentiation; see \cite{MP} and \cite{PBC} for examples. Based on these prior approaches, we present here two methods for reducing noises inherent in the dataset and reconstructing the ground-truth data for later use.

\subsection{Robust smoothing methods}
Several standard smoothing techniques are available for improving regularity of the quantities of interest; e.g. the {\it moving average method}; the LOESS (local regression using weighted linear least squares and a second-degree polynomial model); and the LOWESS (local regression using weighted linear least square and first-degree polynomial model). We refer the reader to  \cite{Cleveland1990} and \cite{Cleveland1988} for more elaborated discussions. However, these smoothing techniques, which are mostly based on least-square regression, are known to be sensitive to outliers, which should ideally be removed before smoothing is applied. A few generic methods for detecting outliers in a dataset are summarized by \cite{Blatna}. For the time series of vehicle location, velocity and acceleration in the NGSIM data, removal of the outliers is done in \cite{MP} and \cite{PBC}.

For the  application described in this paper, we utilize the robust versions of the ``LOESS" or "LOWESS" methods, namely ``RLOESS" or ``RLOWESS", which are implemented and documented in the MATLAB Curve Fitting toolbox  \footnote{URL: http://www.mathworks.com/help/curvefit/smoothing-data.html}. Taking ``RLOESS" as an example, it first smoothes the value at a given point through local regression using nearby data points (the range of those points is determined by a parameter called {\it span}), which is precisely the procedure employed by ``LOESS". Additionally, the ``RLOESS" calculates the residuals from the aforementioned regression and assigns a robust weight to each data point within the span: 
\begin{equation}\label{robustweight}
w_i~=~
\begin{cases}
\left(1-\left( {r_i\over 6 MAD}\right)^2\right)^2\qquad & \hbox{if }~ |r_i|~<~6 MAD
\\
0 \qquad  & \hbox{if }~ |r_i|~\geq~6 MAD
\end{cases}
\end{equation}
where $r_i$ denotes the residual of the $i^{th}$ data point, and MAD ({\it median absolute deviation}) is defined as ${\hbox{median}}(|r_i|)$. The MAD is a measure of how spread out the residuals are. If $r_i$ for a given $i$ is larger than  6MAD, then the $i^{th}$ point is considered an outlier and is assigned weight $0$. Finally, the regression is applied again to smooth the data, but with the robust weights \eqref{robustweight}.

\begin{remark}
If one applies the robust local regression method to location, speed, and acceleration data separately, different degree of smoothing (span) should be used. The reason is that the regularities of these quantities are decreasing, and with each differentiation, the error becomes magnified into increasingly larger ones. Moreover, applying the smoothing separately means that these quantities may no longer be consistent as the derivative of the smoothed location is no longer the smoothed speed, and etc. Such an inconsistency is illustrated further in Figure \ref{figfig1}.
\end{remark}

\subsection{Regularization of numerical differentiation}\label{secrobustdiff}

It is conventionally acknowledged that finite-difference approximations, such as those expressed by \eqref{approxv} and \eqref{approxmv}, tend to amplify noises existing in the subject of numerical differentiation. In some applications, smooth-then-differentiate or differentiate-then-smooth may not yield satisfactory results \citep{Chartrand}. Instead of applying the smoothing method to individual trajectories of location, velocity, or acceleration, one alternative is to utilize a numerical differentiation scheme that is robust against noises in the data. One such robust scheme is Tikhonov regularization \citep{Tikhonov, Chartrand}, and it works as follows. Given a function $f$, its derivative on $[0,\,T]$, denoted by g, is the solution of the following minimization problem. 
\begin{equation}\label{regdiffmin}
\min_{g} ~\mathcal{F}(g)~=~\alpha R(g) + DF (\mathcal{I} g - f)
\end{equation}
where $R(g)$ is a regularization term that penalizes the irregularity in $g$, and it is chosen as the square of the $L^2$ norm of the derivative of $g$ in our study: 
$$
R(g)~\doteq~ \int_0^T \left| g'(\tau)\right|^2\,d\tau
$$
Such a choice forces $g$ to be continuous and reasonably regular by penalizing large values of $|g'|$. The operator $\mathcal{I}$ represents integration, and $\mathcal{I}g (t)\doteq\int_0^t g(\tau)\,d\tau$, $t\in[0,\,T]$. $DF(\mathcal{I} g -f)$ is a {\it data fidelity} term that penalizes the discrepancy between $\mathcal{I}g$ and $f$, and it is chosen as the square of the $L^2$ norm of $\mathcal{I}g- f$.  $\alpha>0$ is a weighting parameter that balances between the two terms.

In a discrete-time setting, denote by $f_i$ and $g_i$ the discrete values of $f$ and $g$ respectively,  $i=1,\ldots, n$. The integration operator is approximated by a rectangular or trapezoidal quadrature and can be written as an affine transformation \footnote{Most numerical integration or differentiation operators are linear in the discrete values of the subject function, therefore they can be represented by affine transformations.} in discrete time:
$$
\Big\{ (\mathcal{I}g)_i:~1\leq i\leq n\Big\}~=~ A g + b \qquad A\in\mathbb{R}^{n\times n},~~ b\in\mathbb{R}^n
$$
where $g=\{g_i:\,1\leq i\leq n\}$. Moreover, the numerical differentiation $g'$ can  be also represented by an affine transformation: 
$$
\Big\{ (g')_i:~1\leq i\leq n\Big\}~=~ C g + d \qquad C\in\mathbb{R}^{n\times n},~~ d\in\mathbb{R}^n
$$
where the dimensions of $A,\,C$ and $b,\,d$ may vary depending on the numerical differentiation scheme.  The discretized minimization problem becomes a quadratic program:
\begin{equation}\label{quadprog}
\min_{g}~\alpha \delta t\cdot (Cg+d)^T (Cg+d) + \delta t \cdot  (Ag+b -f)^T (Ag+b-f)
\end{equation}
with some minor constraints whenever applicable (such as the non-negativity of $v$ when differentiating $x$), where $\delta t$ denotes the time step size. This mathematical program can be solved efficiently with standard solvers given that the dimension of $f$ or $g$ is not too high (no more than a few thousands), which is indeed the case for the NGSIM dataset.

\begin{remark}
Other choices of the regularization term $R(g)$ include the $L^1$ norm of the derivative of $g$, i.e.
\begin{equation}\label{Rchoice1}
R(g)~=~\int_0^T \left| g'(\tau)\right| \,d\tau
\end{equation}
In this case, due to the presence of the absolute value, the discrete optimization problem becomes non-smooth.  \cite{Chartrand} solves this case by applying first-order necessary conditions to the continuous-time problem \eqref{regdiffmin}, which are then discretized. Another choice for $R(g)$ is the square of the $L^2$ norm of $g$ \citep{Tikhonov}:
\begin{equation}\label{Rchoice2}
R(g)~=~\int_0^T \left| g(\tau)\right|^2 \,d\tau
\end{equation}
\noindent which suppresses large values of $g$ without necessarily forcing it to be regular. A quadratic program can be formulated in this case by taking $C$ to be the identity matrix and $d=0$ in \eqref{quadprog}.

\end{remark}

\subsection{Implementation of the denoising methods}
This section presents a numerical study of the two denoising methods mentioned above. We apply the two methods to the NGSIM dataset with the following chosen parameters. For the robust smoothing method, we apply the ``RLOESS" to location (with span $=0.005$), velocity (with span $=0.01$) and acceleration (with span $=0.03$) separately. We note that by doing so, the trajectories are no longer consistent; e.g. the derivative of the smoothed location is not approximated by the smoothed velocity. Regarding the robust differentiation method, we apply the regularization in differentiating location with $\alpha=0.01$ to get velocity, and then differentiate the acquired velocity with $\alpha=0.05$ to get the acceleration. The resulting trajectories are relatively consistent thanks to the data fidelity term in \eqref{regdiffmin}. Notice that in both methods, we have applied different degrees of smoothness (span in the first method and $\alpha$ in the second method) to location, velocity and acceleration due to the natural degradation of regularity as a result of differentiation.

\begin{figure}[h!]
\centering
\includegraphics[width=\textwidth]{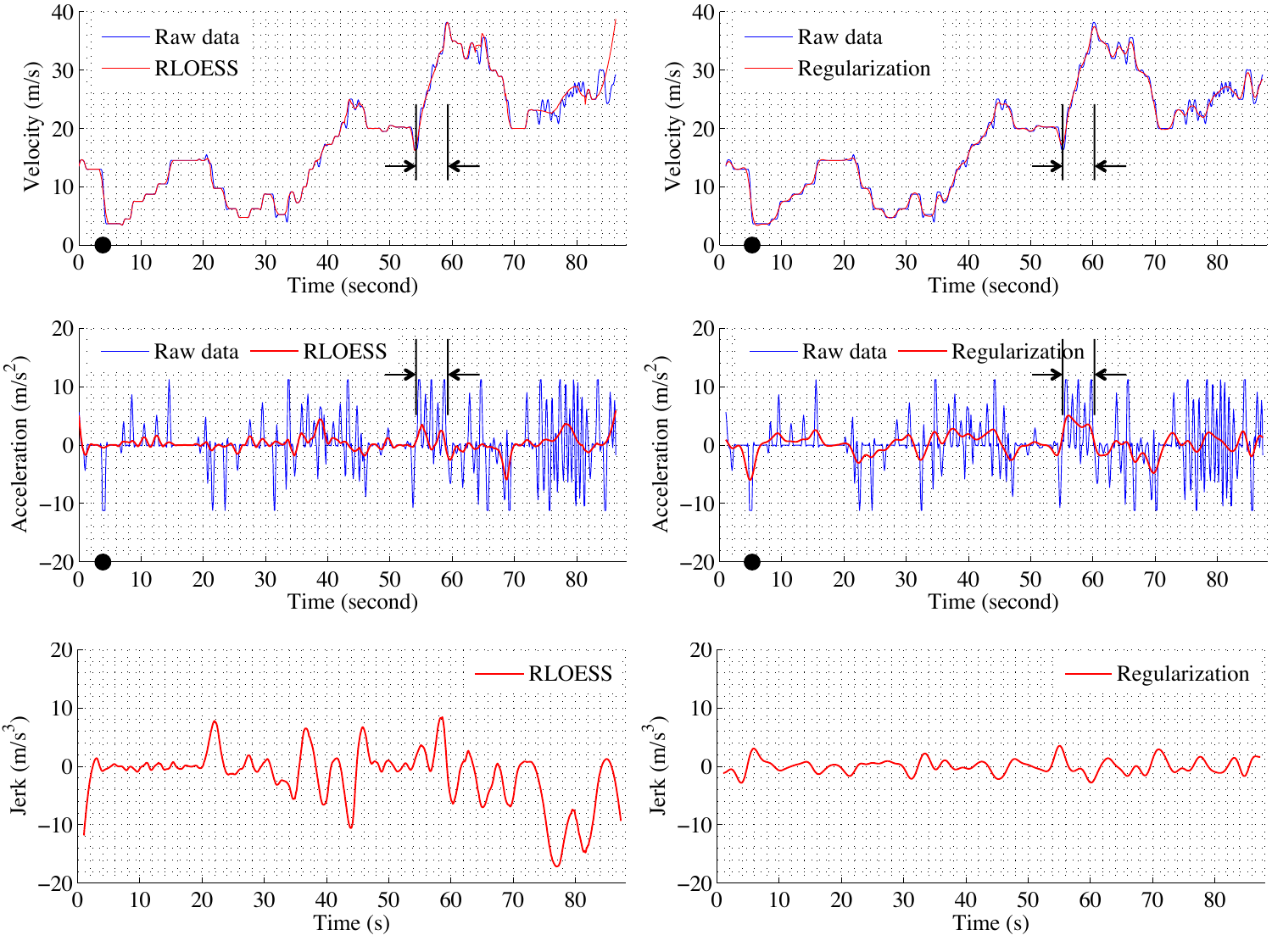}
\caption{Comparison of denoised velocity (top row), acceleration (middle row), and jerk (bottom row) obtained from RLOESS (left column) and regularization (right column). The bottom row only shows the denoised data as the raw data oscillate violently and obstruct the visualization.}
\label{figfig1}
\end{figure}

Figure \ref{figfig1} compares the denoised velocity and acceleration profiles of the same car obtained by the two methods. A discernible difference in the acceleration can be seen while the velocity profiles are similar. We illustrate the qualitative difference between the two acceleration profiles at two time windows. First, at $t\approx 4$s which is marked with a black dot in the figures, the vehicle experienced a sudden velocity reduction, indicating a large negative value in the acceleration. This is, however, insufficiently captured by RLOESS as this single negative value was likely to be treated as an outlier. The regularized differentiation, on the other hand, captures this negative acceleration. The second time window is roughly $[55,\, 60]$ s, during which the car is continuously accelerating. This is reflected by the regularization method, which predicts an acceleration well above zero during this period. However, the RLOESS method yields some negative acceleration values. Overall, it is confirmed that while both methods yield improved regularity in the estimated quantities, the regularization method is more faithful to the physical interpretation of these quantities, as expected from the theory. The RLOESS sometimes treats unusual yet realistic measurements as outliers (such as the sudden speed reduction at $t=4$s). In addition, the regularization method yields acceleration that does not change signs as frequently as the RLOESS, which is consistent with the expected vehicle performance.

In the last row of Figure \ref{figfig1}, we show the time-trajectories of the jerk (rate of change of the acceleration) obtained from the two methods, with span $=0.05$ and $\alpha=0.1$, respectively. An inspection of the overall estimation of velocity, acceleration and jerk reveals that the two methods yield estimation results that are increasingly different when the order of the underlying traffic quantity is higher, which makes sense as the noises are amplified and the regularities are worsened when differentiation is performed. In addition, the regularization method seems to represent the physics of traffic relatively well, while the RLOESS method may be more resilient to the outliers in the data and sensing errors.

For the study presented in the rest of this paper, we re-construct the ground truth using the following two-step method. We first apply the RLOESS method with span $=0.01$ to smooth the location data. Then, we use the smoothed location to perform robust numerical differentiations following Section \ref{secrobustdiff}, and the resulting velocity (with $\alpha=0.01$) and acceleration (with $\alpha=0.03$) are taken as the ground truth. 

\begin{remark}
The proposed two-step method yields location $x$ and acceleration $a$ that are consistent in the sense of regularized differentiation, rather than finite-difference differentiation (i.e. $a_i={x_{i+1}-2x_i+x_{i-1}\over \delta t^2}$). One could devise a further step to construct another location data $\tilde x$ that is consistent with $a$ in the finite-difference sense, although this step is not a necessity for our study. The way to do this is to treat acceleration $a$ as the ground truth and construct $\tilde x$ according to the following optimization problem
$$
\tilde x ~=~\underset{y}{\text{argmin}} \left\| y- x\right\|^2 \qquad\hbox{subject to ~} ~~\mathcal{D}_2y~=~a
$$
where the matrix operation $\mathcal{D}_2 y$ represents the second-order differentiation of $y$ using finite difference \eqref{approxa}, which is expressible as a linear operator. The constraint $\mathcal{D}_2y=a$ is the consistency condition, and the objective is a data fidelity condition that minimizes the discrepancy between $\tilde x$ and the smoothed data $x$. This minimization problem can be formulated as a quadratic program. 

Applying this method will yield a single trajectory $\tilde x$, whose second derivative (in the finite difference sense) is $a$, the smoothed acceleration from our two-step method.

\end{remark}

\section{Model calibration}\label{secNumerical}

Before we implement the computational procedures proposed in Section \ref{secTrafficModel}, we need to estimate the values for $A$, $B$ and $\rho_{jam}$ appearing in \eqref{ptmv}, and subsequent calculations. We use the available dataset to calibrate these modeling parameters.

\subsection{Description of vehicle trajectory data}
\begin{figure}[h!]
\centering
\includegraphics[width=.5\textwidth]{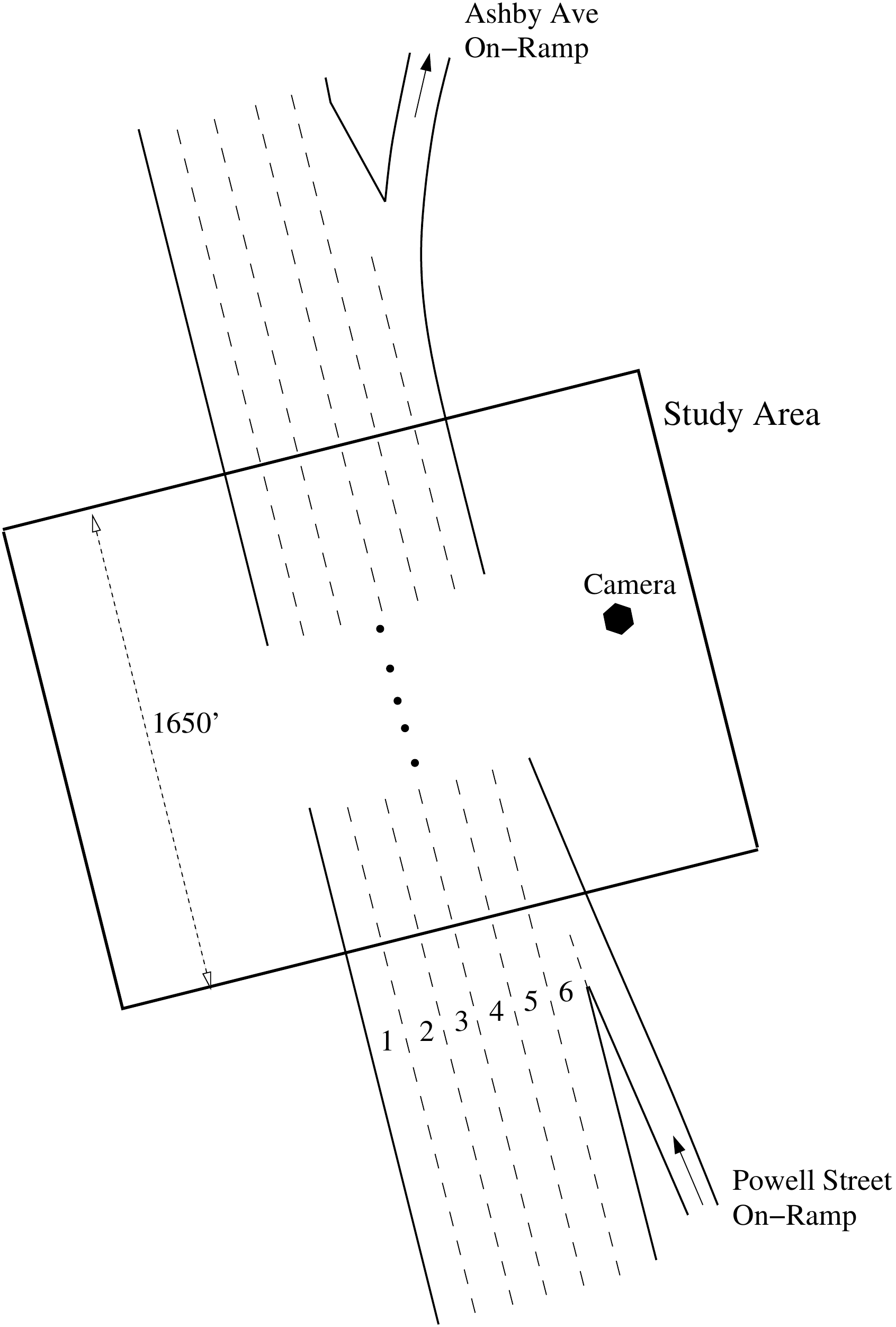}
\caption{The study area, spanning 1650 feet in length in the northbound of Interstate 80 located in Emeryville, CA.}
\label{figStudyArea}
\end{figure}

The NGSIM focuses on the northbound of I-80 located in Emeryville, CA. The highway segment of interest spans 1650 feet in length with an on-ramp at Powell Street and an off-ramp at Ashby Avenue; see Figure \ref{figStudyArea}. The highway segment has six lanes with the leftmost lane being a {\it high-occupancy vehicle} (HOV) lane, and the rightmost one being a merge/diverge lane.  Data were collected using several video cameras. Digital video images were collected over an approximate five-hour period from 2:00 pm to 7:00 pm on April 13, 2005. Complete vehicle trajectories transcribed at a resolution of 1 frame per 0.1 second, along with vehicle type, lane identification and so forth were recorded and processed over three time slots: 4:00 pm - 4:15 pm, 5:00 - 5:15 pm, and 5:15 - 5:30 pm. The layout of the study area is shown in Figure \ref{figStudyArea}.

\subsection{Estimating the density-flow relationship}\label{secmodelcal}
\begin{figure}[h!]
\centering
\includegraphics[width=.35\textwidth]{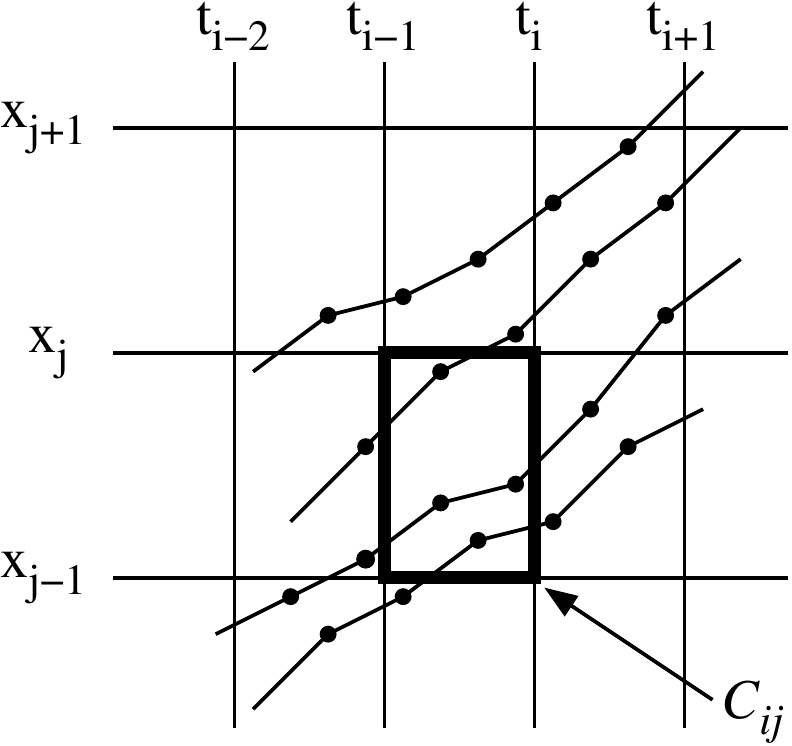}
\caption{Estimation of density and velocity inside a bin. The curves represent vehicle trajectories where the locations are recorded at the solid dots.   For the depicted scenario, the occupancy of bin $C_{ij}$ is three vehicles; the average velocity is taken as the mean of velocities measured/calculated at the four dots inside $C_{ij}$.}
\label{figcells}
\end{figure}

We begin with estimating the density-flow relationship \eqref{funcv2} needed for the congested phase represented by Eqn. \eqref{PTMcongested}. Unlike the LWR model  where the density-flow relation is expressed as a single-valued function, the fundamental diagram corresponding to the congested phase of PTM is a multi-valued map. This  means that a given density $\rho$ corresponds to a continuum range of velocities $v(\rho,\,q),\, q-q^*\in[-1,\,1]$. In order to identify this multi-valued function, we partition the temporal-spatial domain into small bins $C_{ij}\doteq [t_{i-1},\,t_{i}]\times [x_{j-1},\,x_{j}]$, $i =1,\,\ldots,\,N_T$, $j=1,\,\ldots,\,N_X$ where $i$ and $j$ indicate the time step and the spatial step respectively.  The average density associated with $C_{ij}$ is estimated by the number of vehicles whose trajectories indicate their presence in the road segment $[x_{j-1},\,x_{j}]$ during time interval $[t_{i-1},\,t_{i}]$. The velocity inside $C_{ij}$ is calculated as the mean of all velocity measurements collected within this bin. Figure \ref{figcells} illustrates such a procedure. The dimension of the bins used to construct the density-flow relation in the congested phase is  4 (seconds)  $\times$ 400 (feet).

The flow within $C_{ij}$ is calculated as the product of the density and the average velocity. The resulting density-flow and density-velocity plots for the entire study period are shown in Figure \ref{figfds}.

\subsection{Constructing the congested region in the fundamental diagram}\label{secfdcal}

\begin{figure}[h!]
\begin{minipage}[b]{.49\textwidth}
\centering
\includegraphics[width=1\textwidth]{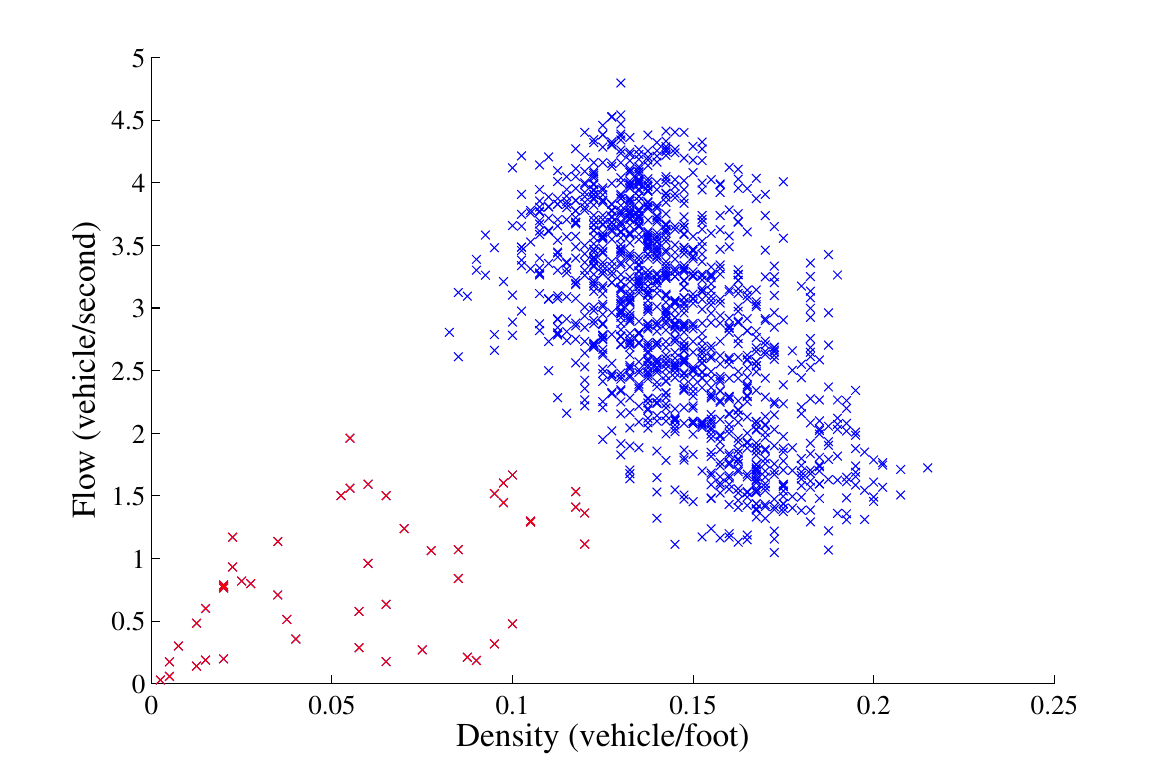}
\end{minipage}
\begin{minipage}[b]{.49\textwidth}
\centering
\includegraphics[width=1\textwidth]{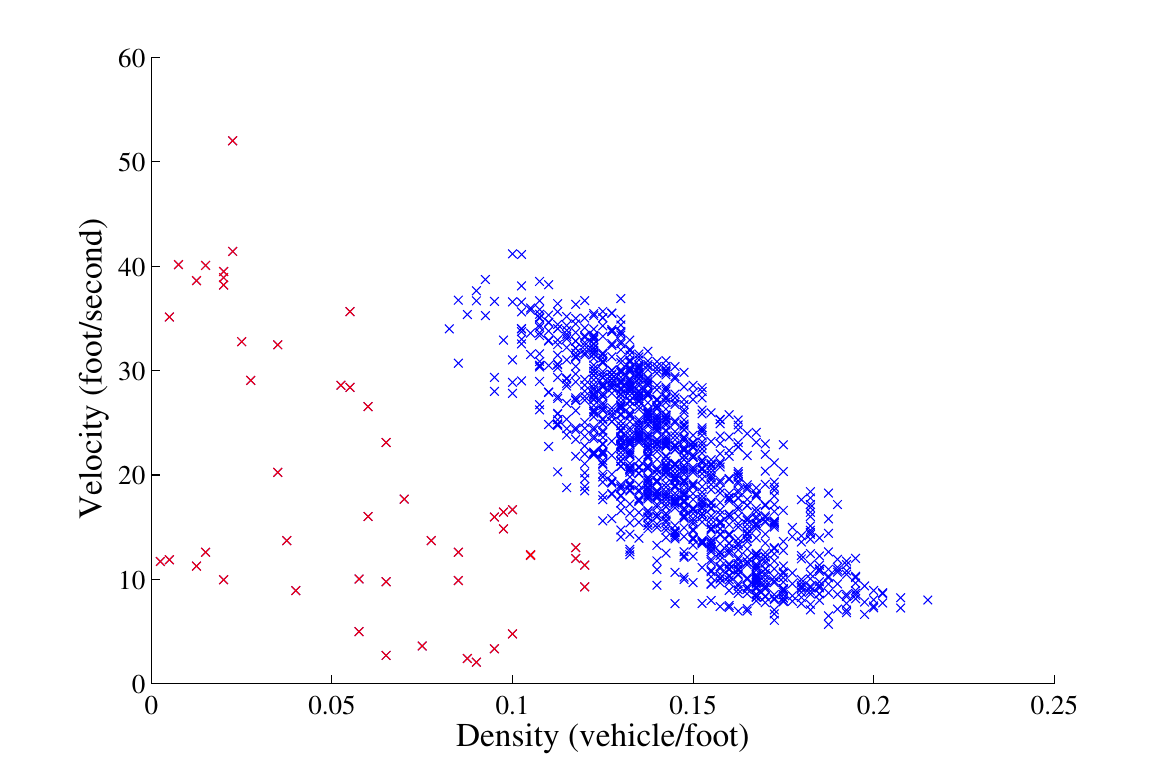}
\end{minipage}
\caption{The congested branch of the density-flow relationship (left) and the density-velocity relationship (right) for the PTM, both expressed as set-valued functions of density. The red crosses represent outliers that are discarded in the curve fitting.}
\label{figfds}
\end{figure}

\begin{figure}[h!]
\begin{minipage}[b]{.49\textwidth}
\centering
\includegraphics[width=1\textwidth]{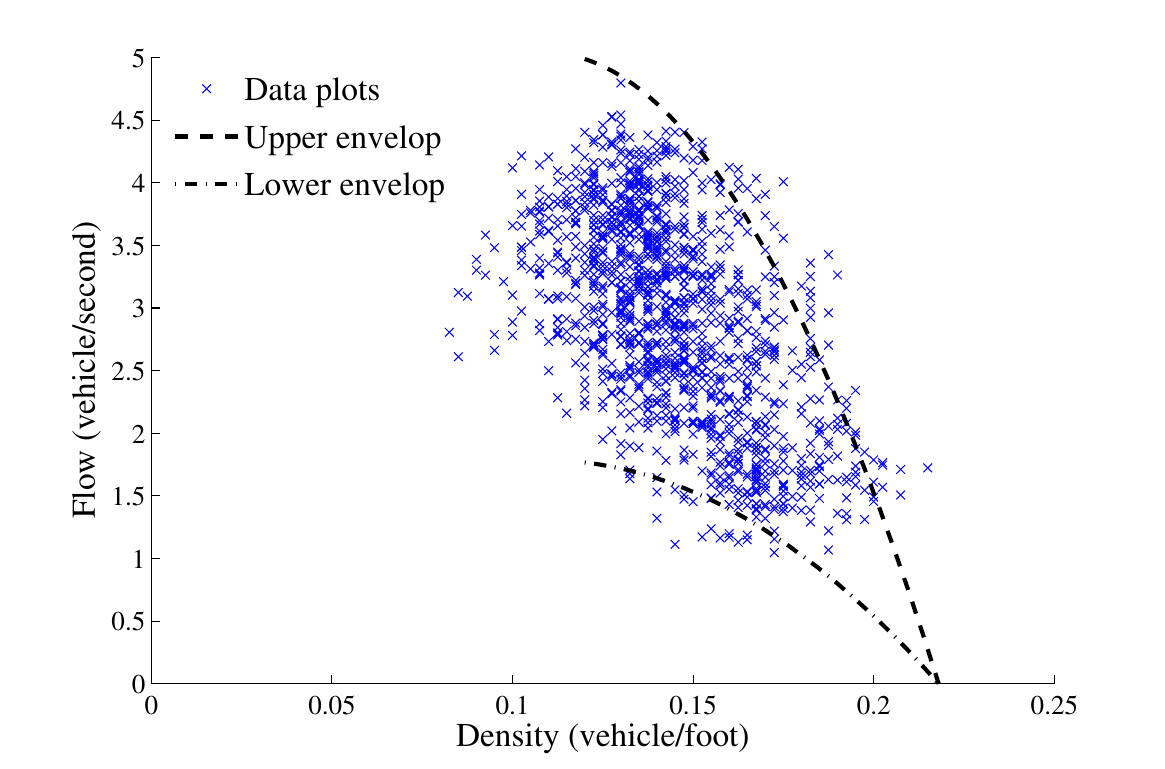}
\end{minipage}
\begin{minipage}[b]{.49\textwidth}
\centering
\includegraphics[width=1\textwidth]{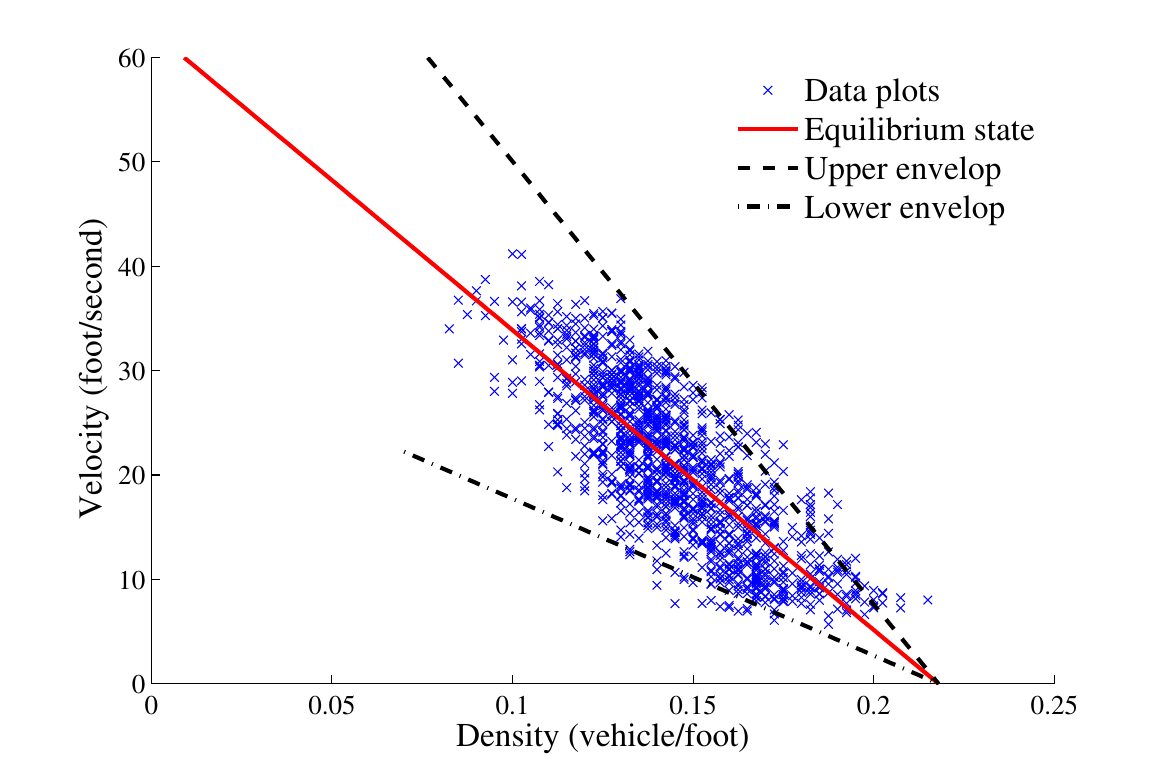}
\end{minipage}
\caption{The fitted density-flow relationship (left) and density-velocity relationship (right) with $A=287.25$ ($\hbox{foot}^2/\hbox{vehicle}/ \hbox{second}$), $B=137$ ($\hbox{foot}^2/\hbox{vehicle}/ \hbox{second}$),  and $\rho_{jam}=0.218$ (vehicle/foot).}
\label{figfdsall}
\end{figure}

Equation \eqref{ptmv} suggests an affine density-velocity relationship $v=A(\rho_{jam}-\rho)$ when the perturbation $q-q^*$ is zero, a situation referred to as the {\it equilibrium state}.  To determine $A$ and $\rho_{jam}$, we conduct the following procedure. First, we visually remove the outliers in the density-velocity plot; see Figure \ref{figfds}. Then, we conduct a  linear regression to determine the equilibrium density-velocity relationship: $v=A(\rho_{jam}-\rho)$. The best fit is $v=-287.25 \rho+ 62.61$, which implies $A=287.25$ ($\hbox{foot}^2/\hbox{vehicle}/ \hbox{second}$) and $\rho_{jam}=0.218$ (vehicle/foot). The linear fit is indicated as the red solid line in the right part of Figure \ref{figfdsall}.

It remains to determine $B$, which effectively dictates the width of the 2-D area in the congested domain of the fundamental diagram.  Since $q-q^*\in[-1,\,1]$, the upper and lower envelops of the congested domain in the density-velocity relationship, depicted in the right part of Figure \ref{figfdsall}, are respectively
$$
\begin{cases}A(\rho_{jam}-\rho)+B(\rho_{jam}-\rho)\qquad &\hbox{(upper envelop)} \\
A(\rho_{jam}-\rho)-B(\rho_{jam}-\rho)\qquad &\hbox{(lower envelop)} 
\end{cases}
$$
Since $A$ and $\rho_{jam}$ are known, we choose the smallest value of $B$ such that $>95\%$ of the data points fall within the area formed by the upper and lower envelopes; see the right part of Figure \ref{figfdsall}. This yields $B=137$ ($\hbox{foot}^2/\hbox{vehicle}/ \hbox{second}$), and the resulting envelops are shown in the right part of Figure \ref{figfdsall}.
Consequently, the upper and lower envelops  in the density-flow relationship depicted in the left part of Figure \ref{figfdsall} are uniquely determined as
$$
\begin{cases}A(\rho_{jam}-\rho)\rho+B(\rho_{jam}-\rho)\rho \qquad &\hbox{(upper envelop)} \\
A(\rho_{jam}-\rho)\rho-B(\rho_{jam}-\rho)\rho \qquad &\hbox{(lower envelop)} 
\end{cases},
$$
which are quadratic curves (parabolas); see the left part of Figure \ref{figfdsall}.

\section{Estimating traffic quantities along vehicle trajectories}\label{secestimatetraffic}
In this section, we present the estimation results associated with various first- and higher-order  traffic quantities along vehicle trajectories. Those quantities include: velocity and acceleration given by \eqref{approxmv}-\eqref{approxa}, vehicle density given by \eqref{eqn3}, \eqref{eqn9} or Algorithm \ref{alg1},  deviation given by \eqref{eqn4}, \eqref{eqn10} or Algorithm \ref{alg1},  and the power demand and HC emission rate given by  \eqref{ztot}-\eqref{hcrate}. Different sampling frequencies will be used to investigate the deterioration of estimation quality caused by under sampling. Recall that the vehicle locations in the NGSIM raw dataset are recorded every 0.1 second. In order to accommodate different sampling frequencies, we consider an integer $N$, and then extract data from the raw dataset every $N$ points. For example, $N=30$ implies a sampling period of $0.1 \times 30=3$ seconds. We then compute the relative $L^1$ error between the quantities estimated with $N>1$ and the ground truth.

All the numerical results reported below are based on data collected during 4:00 pm - 4:15 pm and 5:00 pm - 5:30 pm. The total numbers of vehicles involved in these time periods is 5677. The results summarized in tables or figures below are based on data in the entire study period of 45 minutes unless otherwise stated. As we have explained in detail in Section \ref{secsmooth},  the ground-truth quantities used throughout our numerical experiments are re-constructed using a combined approach of smoothing and robust differentiation.

\subsection{Velocity and acceleration}\label{secNumericalva}

\begin{figure}[h!]
\centering
\includegraphics[width=\textwidth]{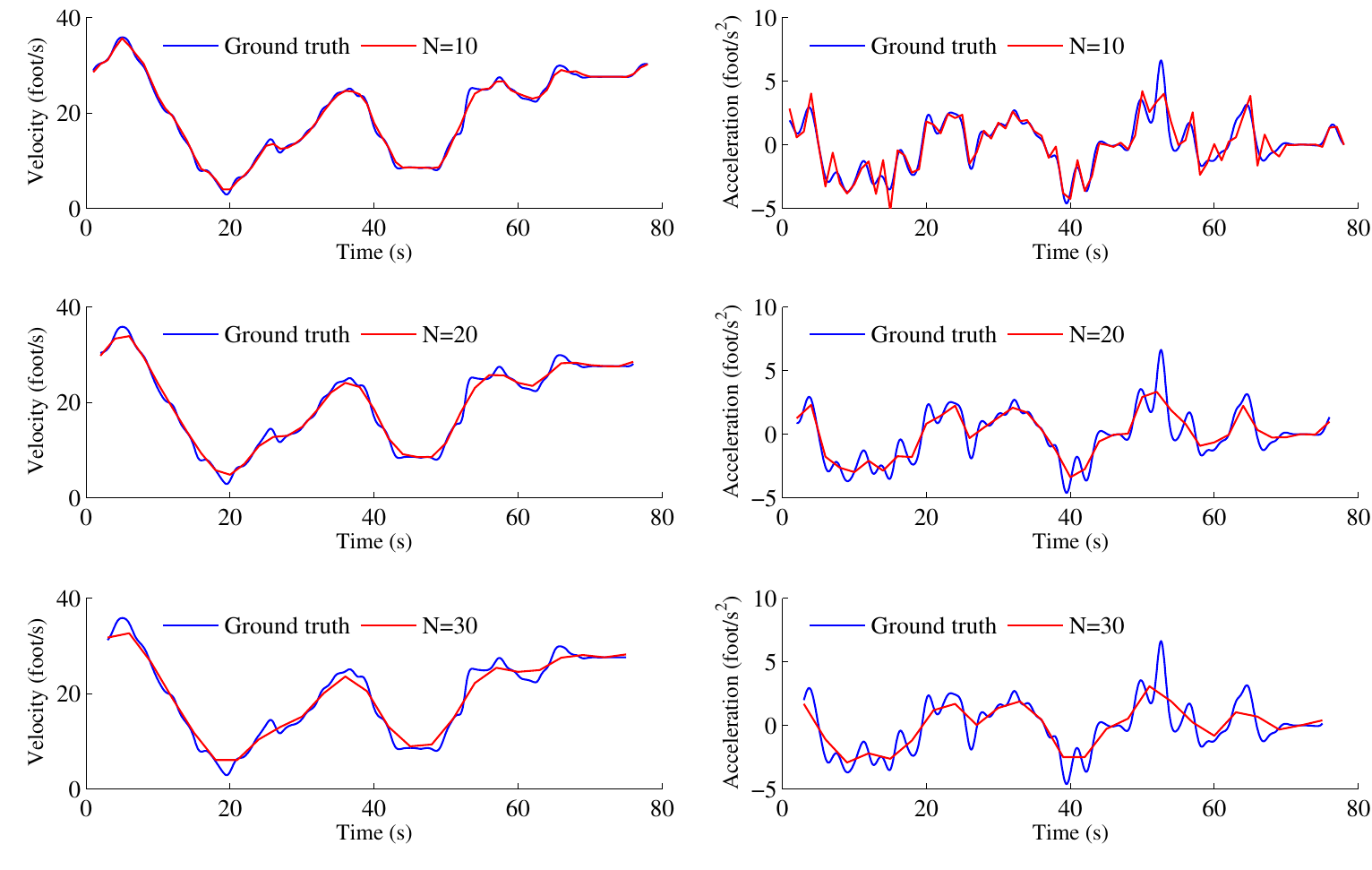}
\caption{The ground-truth vs. estimated velocity (left column) and acceleration (right column) of the same vehicle, using different sampling frequencies ($N=10,\,20,\,30$).}
\label{figacceleration}
\end{figure}

Vehicle speed and acceleration are among the most fundamental traffic quantities. Following the steps explained at the beginning of this section, we compute and summarize in Table \ref{tabva} the mean and standard deviation of the relative $L^1$ errors for the velocity and acceleration. Overall, the estimation of both these quantities deteriorates with under sampling. Moreover, the acceleration estimation is more susceptible to under sampling than velocity estimation. This is because much of the higher-order variation in the acceleration takes place on a smaller time scale, as confirmed by the poor regularity in the acceleration profile (see, for example, Figure \ref{figfig1} and Figure \ref{figacceleration}), and cannot be sufficiently captured with under sampling. Figure \ref{figacceleration} shows one example of the the reconstructed velocity and acceleration of the same car by using different values of $N$. It is confirmed that the fitting of velocity  shows an overall improvement over the fitting of acceleration, which is caused by the poorer regularity in the latter quantity, i.e. more oscillations are observed in a short time period.

\begin{table}[h!]
\begin{center}
\begin{tabular}{|c|c|c|c||c|c|c|c|}\hline
& \multicolumn{3}{|c||}{ Mean (\%)}& \multicolumn{3}{|c|}{Standard deviation (\%)}\\
\hline
$N$   & 10 & 20 & 30 & 10 & 20 & 30\\\hline
$\rule{0pt}{20pt} \displaystyle{\big\|  v_{true}-v_{N}\big\|_{L^1}\over \big\| v_{true}\big\|_{L^1}}$  & 2.15 & 4.35  & 6.41                  & 0.63 &  1.46  & 2.53   \\
\hline
$\displaystyle{\big\| a_{true}-a_{N}\big\|_{L^1}\over \big\| a_{true}\big\|_{L^1}}$ &
 20.91 & 47.43  & 63.51                   & 6.23 & 11.49  & 11.94 \\\hline
\end{tabular}
\end{center}
\caption{Relative errors of velocity and acceleration estimation using different sampling frequencies. $v_{true}$ and $a_{true}$ denote  the ground-truth velocity and acceleration respectively; $v_{N}$ and $a_{N}$ denote the reconstructed velocity and acceleration with different values of $N$.}
\label{tabva}
\end{table}

\subsection{Vehicle density and deviation}\label{secNumericalrhoq}

The density $\hat \rho$ and deviation $\hat q$ are estimated by \eqref{eqn3}-\eqref{eqn4} (Method 1),  \eqref{eqn9} -\eqref{eqn10} (Method 2), or  Algorithm \ref{alg1} (Method 3) under various assumptions on the source term and traffic flow conditions. The means and standard deviations of the relative $L^1$ errors are summarized in Table \ref{tabssrhoqentireperiod}.

\begin{table}[h!]
\begin{center}
\small{
\begin{tabular}{|c|c|c|c|c|c|c|c|}\hline
& & \multicolumn{3}{|c|}{ Mean (\%)}& \multicolumn{3}{|c|}{Standard deviation (\%)}
\\
\hline
 & $N$  &  10 & 20 & 30  & 10 & 20 & 30
\\
\hline
\multirow{3}{*}{$\displaystyle{\big\| \hat \rho_{true}-\hat \rho_{N}\big\|_{L^1}\over \big\| \hat \rho_{true}\big\|_{L^1}}$}  
& Method 1  & 2.15 & 4.55  & 6.64            & 0.63 & 1.62 & 2.70 
 \\
 & Method 2    & 3.56 & 6.70  & 9.35          & 1.35 & 2.73  & 4.08
\\
& Method 3   & 7.19 & 13.80  & 19.39          & 7.48 & 14.37  & 19.46
\\
\hline
\multirow{3}{*}{$\displaystyle{\big\| \hat q_{true}-\hat q_{N}\big\|_{L^1}\over \big\| \hat q_{true}\big\|_{L^1}}$}  
& Method 1   & 20.91 & 54.41  &69.32         & 6.23 & 10.96 & 10.23 
 \\
 & Method 2 & -- & --  & --           & -- & --  & --
\\
& Method 3   & 12.80 & 27.58  & 40.75          & 6.30 &  13.43 & 19.22
\\
\hline
\end{tabular}
}
\end{center}
\caption{Estimation of $\hat \rho$ and $\hat q$ based on  the phase transition model with Methods 1--3. $\hat \rho_{true}$ and $\hat q_{true}$ represent the ground-truth quantities estimated using the ground-truth location, velocity and acceleration re-constructed in Section \ref{secsmooth}.  $\hat\rho_{N}$ and $\hat q_{N}$ denote the  estimated quantities with under sampling.  The estimation of $\hat q$ with Method 2 is not shown as the errors are above $100\%$.}
\label{tabssrhoqentireperiod}
\end{table}

All three methods show qualitatively similar results: The estimation of the first-order quantity (density) is more accurate than the estimation of the second-order quantity (deviation) under the same condition. Moreover, both estimates deteriorate as the sampling frequency decreases. These results show some similarities to the estimation of velocity and acceleration (Section \ref{secNumericalva}). Furthermore, method 2 shows slightly worsened yet satisfactory estimation for density $\hat \rho$. On the other hand, its estimation of $\hat q$ suffers significantly with under sampling (all errors are greater than 100\% and thus are not shown in the table). As we have explained in Eqn \eqref{vanishing} that the finite-difference approximation of the measurement ${Dv\over Dt}-vv_x$ is zero. Even with modified finite-difference schemes such as \eqref{modifieddiff} or robust differentiation schemes mentioned in Section \ref{secsmooth}, such a quantity is likely to be very small, causing Eqn \eqref{eqn10} to yield $\hat q$ with small magnitude and small $L^1$ norm, and consequently, leading to large relative $L^1$ errors. This shows that Method 2 is not robust enough in handling the deviation $\hat q$.

Interestingly, the estimation errors for $\hat q$ in the no-source case (Method 3) are much lower than Method 1. This is due to the exact approach employed by Method 3 without simplification regarding the traffic dynamics, as opposed to the other two approaches. However, recall that Method 3 works only when $v>{4\over 9}A\rho_{jam}\approx 27.83$ (foot/second). As shown by the right half of Figure \ref{figfds}, the majority of speed measurements fall below this value. In order to apply this method, we select vehicles whose trajectories indicate at least one non-trivial period (over 20 seconds),  during which the speed is uniformly above $27.83$ foot/second, and then apply our computational procedure (Algorithm 1) to these periods. The total number of vehicles involved in the computation is 1628 (out of 5677 in total). We conclude that the applicability of Algorithm \ref{alg1} for the no-source case is limited, although it yields a more accurate estimation of the deviation than the other two methods.

\begin{figure}[h!]
\centering
\includegraphics[width=\textwidth]{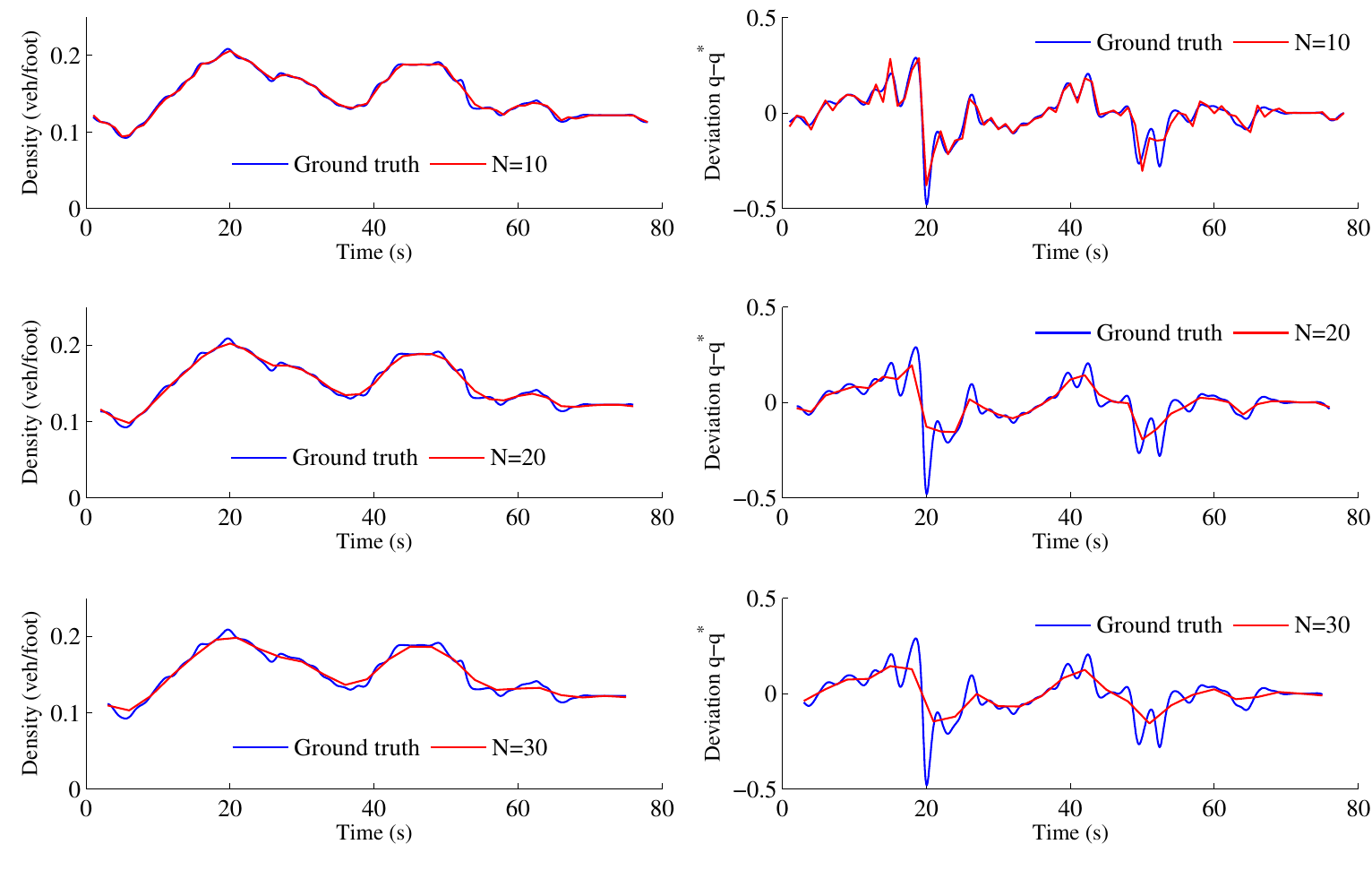}
\caption{The ground-truth and estimated density (left column) and deviation (right column) profiles along the same vehicle trajectory, using different sampling frequencies.}
\label{figRQ}
\end{figure}

Figure \ref{figRQ} shows one example of the reconstructed density and deviation profiles along a single vehicle trajectory, for three different sampling frequencies ($N=10,\,20,\,30$). First of all, a comparison of these cases leads to a similar conclusion regarding the deterioration of estimating first- and second-order quantities as in Section \ref{secNumericalva}. Moreover, one can clearly observe, from both ground-truth and reconstructed $q-q^*$, that the vehicle frequently experienced changes in traffic conditions, between below-equilibrium speed ($q-q^*<0$) and above-equilibrium speed ($q-q^*>0$). Such a characteristic of congested, unstable traffic cannot be captured by first-order models such as the LWR model as it describes only the equilibrium speed as a function of density.

\subsection{Emission rate and power demand function}\label{secNumericalp}
The HC emission rate and the power demand function are estimated according to \eqref{hcrate} and \eqref{ztot}  respectively, where the velocity $v$ and acceleration $a$ are approximated by \eqref{approxmv} and \eqref{approxa}.  Table \ref{tabpowerdemandundersampling} summarizes the relative $L^1$ errors in estimating the HC emission rate $r^{HC}$ and the power demand $Z$, when different sampling frequencies are employed.  

\begin{figure}[h!]
\centering
\includegraphics[width=\textwidth]{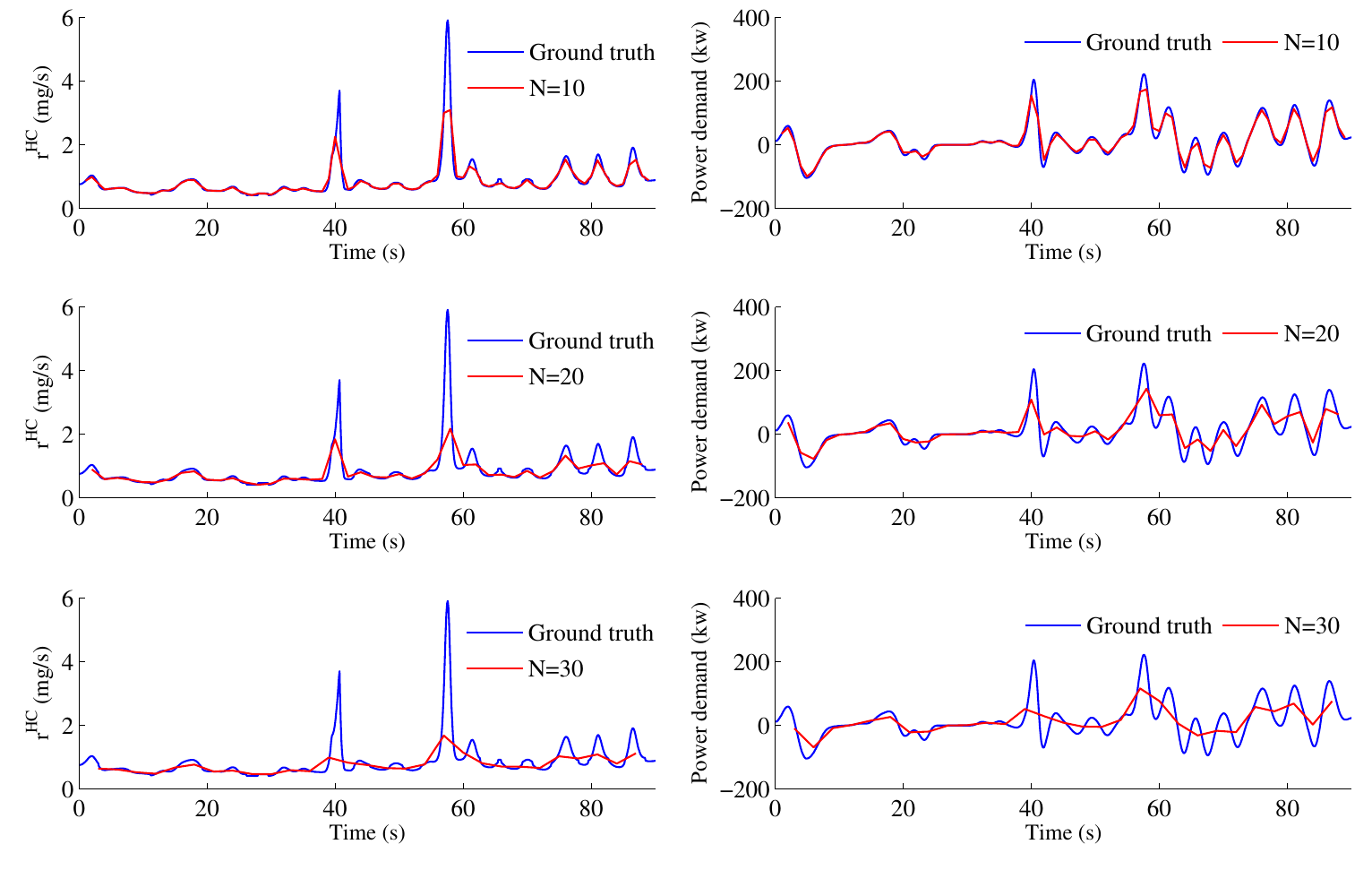}
\caption{Ground-truth and estimated hydrocarbon emission rate (left column) and power demand function (right column).}
\label{figEP}
\end{figure}

Overall, estimation of these higher-order quantities suffer from under sampling since higher-order variations in the acceleration profiles are ignored, as expected from Section \ref{secNumericalva}. Moreover, the deterioration in the accuracy seems most significant when the sampling period increases from $N=10$ (1 second) to $N=20$ (2 seconds), which suggests that the sampling period should ideally be well below 2 seconds in order to capture these quantities. This can be confirmed from Figure \ref{figEP}, which shows both ground-truth and estimated $r^{HC}$ and $Z$ associated with the same vehicle. We see that while $N=10$ captures the variations in emission rate and power demand relatively well, $N=20$ or $30$ causes the estimation to miss the majority of the peak values. Finally, a comparison between $r^{HC}$ and $Z$ for the same value of $N$ reveals that the errors in the former case is much smaller. The reason is that while both estimates misinterpret the high variations by similar amount, the time-trajectory of $r^{HC}$ is uniformly above some positive value while the time-trajectory of $Z$ remains centered at zero. Therefore, the relative $L^1$ error for $r^{HC}$ tends to be smaller.

\begin{table}[h!]
\begin{center}
\small{
\begin{tabular}{|c|c|c|c|c|c|c|}\hline
 & \multicolumn{3}{|c|}{ Mean (\%)}& \multicolumn{3}{|c|}{Standard deviation (\%)}
\\
\hline
  $N$   & 10 & 20 & 30  & 10 & 20 & 30
\\
\hline
   $\displaystyle{\big\| r^{HC}_{true}-r^{HC}_{N}\big\|_{L^1}\over \big\|r^{HC}_{true}\big\|_{L^1}}$  & 13.39 & 23.58  & 27.30            & 16.77 & 24.69 & 24.70
 \\
\hline

  $\displaystyle{\big\| Z_{true}-Z_{N}\big\|_{L^1}\over \big\|Z_{true}\big\|_{L^1}}$     & 18.69 &  44.93  & 60.42            & 4.95 & 9.36  & 10.47
 \\
\hline
\end{tabular}}
\end{center}
\caption{Estimations of the power demand functions $Z$, the hydrocarbon emission rate $r^{HC}$, and the fuel consumption rate $r^{FC}$ based on vehicle trajectories.  Subscripts $true$ and $N$ denote the ground-truth quantity and the reconstructed quantity, respectively.}
\label{tabpowerdemandundersampling}
\end{table}

\section{Estimating Eulerian quantities}\label{secEulerian}
The previous section is mainly concerned with traffic quantities associated with a Lagrangian particle; i.e. a moving car.  It would be desirable to further explore the effect of under sampling in an Eulerian framework; that is, we will study traffic quantities represented on the spatial-temporal domain. We are also prompted to examine how the estimating error depends on the penetration rate of probe vehicles (mobile sensors), in addition to the sampling frequency. The Eulerian traffic quantities studied in this section include vehicle density and deviation. Vehicle emission rate will be studied in the next Section.

\subsection{Estimating vehicle densities}\label{seccellrho}

The procedure for estimating Eulerian density based on vehicle trajectory data and the phase transition model is illustrated as follows. We consider again the temporal-spatial bins $C_{ij},\, i=1,\ldots, N_T$, $j=1,\ldots, N_X$, each expressed as a product of intervals $[t_{i-1},\,t_{i}]\times [ x_{j-1},\, x_{j}]$.  Recall that given a discrete-time trajectory of a vehicle:
$$
\ldots,\, x(\tau_{k-1}),\,x(\tau_k),\,x(\tau_{k+1}),\,\ldots
$$
where $\{\tau_k\}$ is a fixed time grid, one can estimate the velocity $v(\tau_k)$, acceleration $a(\tau_k)$, and the Lagrangian density $\rho\big(\tau_k,\,x(\tau_k)\big)$ using techniques elaborated in Section \ref{secMeasurement} and Section \ref{secTrafficModel}. In order to estimate the average density of $C_{ij}$, we search for the probe vehicles whose trajectories intersect $C_{ij}$. Then the quantity $\rho_{ij}$ associated with $C_{ij}$ is estimated as the average of $\rho(\tau_k,\,x(\tau_k))$ for all $k$ such that $(\tau_k,\,x(\tau_k))\in [t_{i-1},\,t_{i}]\times[x_{j-1},\,x_{j}]$.

\begin{figure}[h!]
\centering
\includegraphics[width=.8\textwidth]{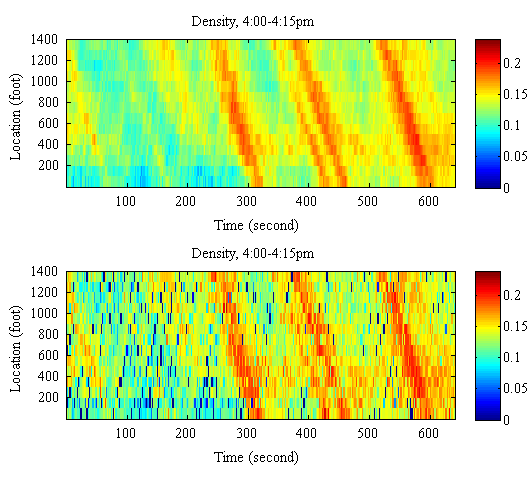}
\caption{Reconstruction of cell densities (in vehicle/foot) during 4:00 pm - 4:15 pm. The top figure shows the density obtained with $N=1$ and 100\% penetration rate. The bottom figure  shows the estimated density with $N=30$ and penetration rate $=50\%$. The inactive bins are assigned zero density and marked as dark blue in the figure.}
\label{figdensity_combined}
\end{figure}

One example of such a calculation is presented in Figure \ref{figdensity_combined}, where we utilize  vehicle trajectories with a sampling period of $N=30$ (3 seconds) and a 50\% probe penetration rate to perform the estimation of bin-based densities, which is then compared with the densities obtained from $N=1$ and a 100\% penetration rate.  Figure \ref{figdensity_combined} shows consistency between the two density estimations. Notice from the bottom of Figure \ref{figdensity_combined} that a few bins are marked with dark blue, which means that they do not have any valid location measurements that support the computation of density values.  We call these bins {\it inactive}, which are automatically assigned zero density and indicated by dark blue in Figure \ref{figdensity_combined}. It is expected that as the penetration rate and the sampling frequency decrease further, the number of inactive bins will increase drastically. To avoid this, one should choose the dimension of a bin reasonably large to ensure having at least some valid measurements inside each bin. In our numerical results presented in the next section, the dimension of a bin is chosen to be 4 (seconds) $\times$ 400 (feet).  We shall use the phase transition model and the LWR model to perform Eulerian density estimation with a range of sampling frequencies and probe penetration rates.

\subsection{Density estimation: A comparative study of the PTM and the LWR mdoel}\label{subsecrhocell}
In this section we evaluate the performance of the aforementioned estimation method in the presence of insufficient data coverage, that is, when the penetration rate and the sampling frequency are low. The purpose is to evaluate and justify the effectiveness of the proposed computational procedure, and to identify certain range of penetration rates and sampling frequencies such that the domain of study is sufficiently covered and the estimation error remains reasonably low.

\begin{figure}[h!]
\begin{minipage}[b]{.49\textwidth}
\centering
\includegraphics[width=1\textwidth]{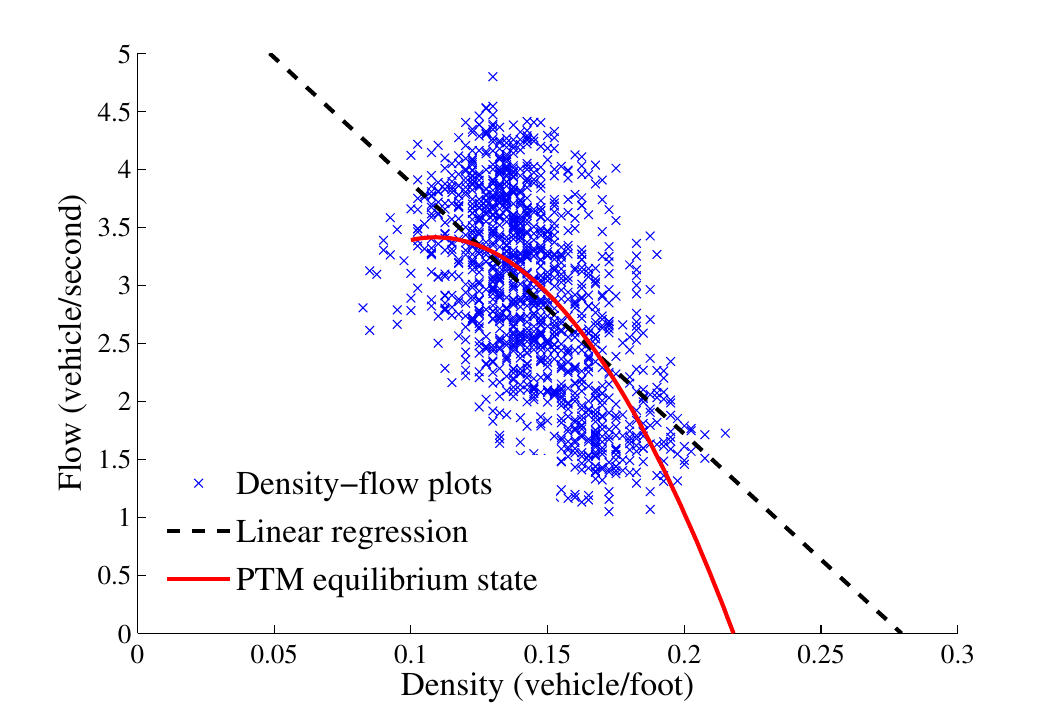}
\end{minipage}
\begin{minipage}[b]{.49\textwidth}
\centering
\includegraphics[width=1\textwidth]{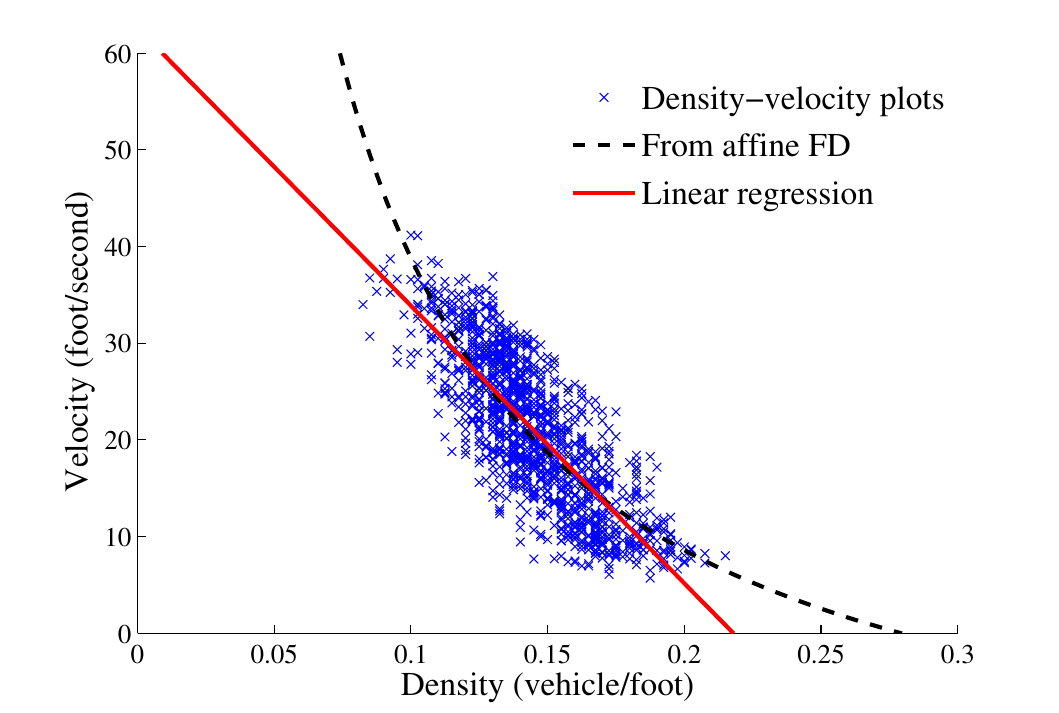}
\end{minipage}
\caption{Fitting of the congested branch of the LWR fundamental diagram. Our study uses the equilibrium density-velocity or density-flow relationship from the PTM case, which are indicated as solid red lines.}
\label{figfdsall_LWR}
\end{figure}

In addition, we propose a comparative study of second-order models (the phase transition model) and first-order models (the LWR model) in terms of their accuracies in density estimation. For the  LWR model \eqref{PTMfree}, we need to calibrate the single-valued fundamental diagram (FD) in view of the  density-flow relationship shown in Figure \ref{figfds}. To do this, we first consider an affine congested branch by applying a least-square linear regression to the density-flow plot after removing the outliers; see the left figure of Figure \ref{figfdsall_LWR}. The result, however,  was rather peculiar (with very large $\rho_{jam}$) as shown in Figure \ref{figfdsall_LWR}. Applying this affine FD leads to substantial error in the estimation and thus is not adopted here. We instead employ the equilibrium density-flow relationship from the PTM case (notice that this relationship is obtained from a least-square linear regression of the density-velocity plots; see the right part of Figure \ref{figfdsall_LWR}): 
$$
\hbox{flow}~=~A\rho(\rho_{jam}-\rho)
$$
where $A=287.25$ ($\hbox{foot}^2$/vehicle/second), $\rho_{jam}=0.218$ (vehicle/foot); see Figure \ref{figfdsall_LWR} for the curve fitting result, where the utilized curves are indicated with solid red lines.

It can be seen that the single-valued density-flow relationship, no matter how derived, does not capture deviation from the equilibrium states, i.e. the scattered measurements in a 2-D region. Therefore, the LWR model treats the traffic stream as stable and homogeneous; as a result, second-order variations such as acceleration/deceleration are not captured.  In order to apply the LWR model to perform density estimation, we  first utilize the probe trajectory data to obtain information of the instantaneous speed $v$; then we evaluate the function $\rho=g^{-1}(v)$ to get the estimated density value, where $g(\rho)$ is the invertible density-speed relationship depicted on the right half of Figure \ref{figfdsall_LWR}. The rest of the procedure is the same as the PTM case.

The ground-truth density associated with each spatial-temporal bin is calculated by counting the number of vehicles present in that bin, using ground-truth trajectory data. The ground-truth density is then compared with the estimated densities obtained with a combination of lower sampling frequencies and lower penetration rates. The relative errors for all active bins are averaged and  shown in Table \ref{tabPTMvsLWR1} for period 4:00-4:15 pm, in Table \ref{tabPTMvsLWR2} for period 5:00-5:15 pm, and in Table \ref{tabPTMvsLWR3} for period 5:15-5:30 pm. Here in these calculations, Method 1 is employed by the PTM-based estimation. The results with Method 2 is qualitatively similar and is omitted here.

From all three tables, we notice that the penetration rate of probe vehicles has a stronger effect on the estimation accuracy than the sampling frequency. This again confirms our observation that, when it comes to estimating first-order quantities, under sampling does not affect the accuracy as much as the penetration rate.  On the other hand, the substantial variation in the penetration rate, from $100\%$ to $2\%$, has caused approximately 14\%  difference in both PTM-based and LWR-based estimations. Moreover, the PTM-based estimation is more accurate than the LWR-based estimation for any given combination of penetration rate and sampling frequency, and the difference is more evident for lower penetration rates. This is explained by the fact that first-order models only captures the equilibrium state of the density-flow or density-velocity relationship, and ignores possible deviations from such an equilibrium state. As shown in Figure \ref{figfdsall_LWR}, the congested traffic is unstable and likely to deviate from the  equilibrium state. Thus, the second-order model (PTM) is able to capture these deviations through the use of the variable $q$ and a multi-valued fundamental diagram. As a result, the second-order model is able to estimate vehicle densities more accurately than the LWR model for congested highway traffic.

\begin{table}[h!]
\begin{center}
\begin{tabular}{ccc|c|c|c|c|c|c|}  \cline{4-9}
 & & & \multicolumn{6}{|c|}{Probe Vehicle Penetration Rate}
\\
\cline{3-9}
  &  & \multicolumn{1}{|c|}{\multirow{1}{*}{Sampling}} & \multicolumn{1}{|c|}{\multirow{2}{*}{100 \%}} & \multicolumn{1}{|c|}{\multirow{2}{*}{50 \%}} & \multicolumn{1}{|c|}{\multirow{2}{*}{20 \%}}  &  \multicolumn{1}{|c|}{\multirow{2}{*}{10 \%}} & \multicolumn{1}{|c|}{\multirow{2}{*}{5 \%}} & \multicolumn{1}{|c|}{\multirow{2}{*}{2 \%}}
\\
& & \multicolumn{1}{|c|}{\multirow{1}{*}{Period}}  & \multicolumn{1}{|c|}{\multirow{2}{*}{}} & \multicolumn{1}{|c|}{\multirow{2}{*}{}} & \multicolumn{1}{|c|}{\multirow{2}{*}{}}  &  \multicolumn{1}{|c|}{\multirow{2}{*}{}} & \multicolumn{1}{|c|}{\multirow{2}{*}{}} & \multicolumn{1}{|c|}{\multirow{2}{*}{}} \\

\hline
\multicolumn{1}{|c|}{\multirow{4}{*}{PTM}} &  \multicolumn{1}{c|}{\multirow{2}{*}{Average}}       &     $N=10$   &  8.74     & 9.25     & 11.04    &  13.40    & 16.95 &  21.53   \\
\multicolumn{1}{|c}{\multirow{1}{*}{}} & \multicolumn{1}{|c|}{\multirow{2}{*}{Error (\%)}}                  &  $N=20$   &  8.86   &  9.44    & 11.19    & 13.48     & 17.05  &   21.26 \\                     
\multicolumn{1}{|c}{\multirow{1}{*}{}} & \multicolumn{1}{|c|}{\multirow{1}{*}{}}                                  &  $N=30$    &  8.99    &  9.57   &   11.41   & 13.74    & 17.25   &   21.04\\   
 \cline{2-9}
\multicolumn{1}{|c}{\multirow{1}{*}{}} &   \multicolumn{2}{|c|}{Coverage Rate (\%)}                                            & 100.00   & 100.00   & 100.00   & 99.83   & 97.67  &  79.70        \\
 \hline                                
 
 \multicolumn{1}{|c|}{\multirow{4}{*}{LWR}} &  \multicolumn{1}{c|}{\multirow{2}{*}{Average}} &          $N=10$  &  12.66     &  13.29      & 15.47     &  18.37   & 22.45  &  27.52\\
\multicolumn{1}{|c}{\multirow{1}{*}{}} & \multicolumn{1}{|c|}{\multirow{2}{*}{Error (\%)}}                 &   $N=20$  &  12.58     &   13.24     & 15.41    &  18.29  &  22.33   &  27.04\\                     
\multicolumn{1}{|c}{\multirow{1}{*}{}} & \multicolumn{1}{|c|}{\multirow{1}{*}{}}                                  &   $N=30$   &  12.43    &   13.03      & 15.31    &  18.26  &  22.26  &   26.44\\    
\cline{2-9}
\multicolumn{1}{|c}{\multirow{1}{*}{}} &   \multicolumn{2}{|c|}{Coverage Rate (\%)}                                 & 100.00   & 100.00   & 100.00   & 99.83   & 97.67  &  79.70  \\
 \hline                                                                                                                                   
\end{tabular}
\end{center}
\caption{Comparison of the PTM and LWR models in estimating bin-based densities  using different sampling periods and penetration rates, for the period 4:00-4:15 pm.  ``Coverage Rate" is the ratio of the number of active bins and the total number of bins.}
\label{tabPTMvsLWR1}
\end{table}

\begin{table}[h!]
\begin{center}
\begin{tabular}{ccc|c|c|c|c|c|c|}  \cline{4-9}
 & & & \multicolumn{6}{|c|}{Probe Vehicle Penetration Rate}
\\
\cline{3-9}
  &  & \multicolumn{1}{|c|}{\multirow{1}{*}{Sampling}} & \multicolumn{1}{|c|}{\multirow{2}{*}{100 \%}} & \multicolumn{1}{|c|}{\multirow{2}{*}{50 \%}} & \multicolumn{1}{|c|}{\multirow{2}{*}{20 \%}}  &  \multicolumn{1}{|c|}{\multirow{2}{*}{10 \%}} & \multicolumn{1}{|c|}{\multirow{2}{*}{5 \%}} & \multicolumn{1}{|c|}{\multirow{2}{*}{2 \%}}
\\
& & \multicolumn{1}{|c|}{\multirow{1}{*}{Period}}  & \multicolumn{1}{|c|}{\multirow{2}{*}{}} & \multicolumn{1}{|c|}{\multirow{2}{*}{}} & \multicolumn{1}{|c|}{\multirow{2}{*}{}}  &  \multicolumn{1}{|c|}{\multirow{2}{*}{}} & \multicolumn{1}{|c|}{\multirow{2}{*}{}} & \multicolumn{1}{|c|}{\multirow{2}{*}{}} \\

\hline
\multicolumn{1}{|c|}{\multirow{4}{*}{PTM}} &  \multicolumn{1}{c|}{\multirow{2}{*}{Average}}       &     $N=10$   &  8.38     & 8.78     & 9.60    &  11.27    & 14.19 &  17.92   \\
\multicolumn{1}{|c}{\multirow{1}{*}{}} & \multicolumn{1}{|c|}{\multirow{2}{*}{Error (\%)}}                  &  $N=20$   &  8.49   &  8.90    &  9.71    & 11.38     & 14.21  &   17.83 \\                     
\multicolumn{1}{|c}{\multirow{1}{*}{}} & \multicolumn{1}{|c|}{\multirow{1}{*}{}}                                  &  $N=30$    &  8.85    &  9.23   &  10.01   & 11.63    & 14.47   &   17.76\\   
 \cline{2-9}
\multicolumn{1}{|c}{\multirow{1}{*}{}} &   \multicolumn{2}{|c|}{Coverage Rate (\%)}                                            & 100.00   & 100.00   & 100.00 & 99.97   & 98.57   & 82.88        \\
 \hline                                

 \multicolumn{1}{|c|}{\multirow{4}{*}{LWR}} &  \multicolumn{1}{c|}{\multirow{2}{*}{Average}} &          $N=10$  &  10.78     &  11.26      & 12.24    &  14.24   & 17.85  &  23.09\\
\multicolumn{1}{|c}{\multirow{1}{*}{}} & \multicolumn{1}{|c|}{\multirow{2}{*}{Error (\%)}}                 &   $N=20$  &  10.96     &   11.44     & 12.44    &  14.44  &  17.91   & 22.95\\                     
\multicolumn{1}{|c}{\multirow{1}{*}{}} & \multicolumn{1}{|c|}{\multirow{1}{*}{}}                                  &   $N=30$   &  11.44    &   11.93      & 12.93    &  14.81  &  18.29  &  22.80\\    
\cline{2-9}
\multicolumn{1}{|c}{\multirow{1}{*}{}} &   \multicolumn{2}{|c|}{Coverage Rate (\%)}                                  & 100.00   & 100.00   & 100.00 & 99.97   & 98.57   & 82.88  \\
 \hline                                                                                                                                       
\end{tabular}
\end{center}
\caption{Comparison of the PTM and LWR models in estimating bin-based densities  using different sampling periods and penetration rates, for the period 5:00-5:15 pm. ``Coverage Rate" is the ratio of the number of active bins and the total number of bins.}
\label{tabPTMvsLWR2}
\end{table}

\begin{table}[h!]
\begin{center}
\begin{tabular}{ccc|c|c|c|c|c|c|}  \cline{4-9}
 & & & \multicolumn{6}{|c|}{Probe Vehicle Penetration Rate}
\\
\cline{3-9}
  &  & \multicolumn{1}{|c|}{\multirow{1}{*}{Sampling}} & \multicolumn{1}{|c|}{\multirow{2}{*}{100 \%}} & \multicolumn{1}{|c|}{\multirow{2}{*}{50 \%}} & \multicolumn{1}{|c|}{\multirow{2}{*}{20 \%}}  &  \multicolumn{1}{|c|}{\multirow{2}{*}{10 \%}} & \multicolumn{1}{|c|}{\multirow{2}{*}{5 \%}} & \multicolumn{1}{|c|}{\multirow{2}{*}{2 \%}}
\\
& & \multicolumn{1}{|c|}{\multirow{1}{*}{Period}}  & \multicolumn{1}{|c|}{\multirow{2}{*}{}} & \multicolumn{1}{|c|}{\multirow{2}{*}{}} & \multicolumn{1}{|c|}{\multirow{2}{*}{}}  &  \multicolumn{1}{|c|}{\multirow{2}{*}{}} & \multicolumn{1}{|c|}{\multirow{2}{*}{}} & \multicolumn{1}{|c|}{\multirow{2}{*}{}} \\

\hline
\multicolumn{1}{|c|}{\multirow{4}{*}{PTM}} &  \multicolumn{1}{c|}{\multirow{2}{*}{Average}}       &     $N=10$   &  8.89    & 9.00     & 9.62    &  10.42    & 12.24 &  15.60   \\
\multicolumn{1}{|c}{\multirow{1}{*}{}} & \multicolumn{1}{|c|}{\multirow{2}{*}{Error (\%)}}                &  $N=20$   &   9.02    &  9.13    & 9.74    &  10.53    & 12.35  &  15.47 \\                     
\multicolumn{1}{|c}{\multirow{1}{*}{}} & \multicolumn{1}{|c|}{\multirow{1}{*}{}}                                  &  $N=30$    &  9.17    &  9.29   & 9.93   &   10.71  &   12.52 &   15.39 \\   
 \cline{2-9}
\multicolumn{1}{|c}{\multirow{1}{*}{}} &   \multicolumn{2}{|c|}{Coverage Rate (\%)}                                             & 100.00   & 100.00   & 100.00 & 99.98   & 98.99   & 85.56         \\
 \hline                                
  
 \multicolumn{1}{|c|}{\multirow{4}{*}{LWR}} &  \multicolumn{1}{c|}{\multirow{2}{*}{Average}} &          $N=10$  &  12.69     &  12.87      & 13.49     & 14.45   & 16.82  &  21.03\\
\multicolumn{1}{|c}{\multirow{1}{*}{}} & \multicolumn{1}{|c|}{\multirow{2}{*}{Error (\%)}}                 &   $N=20$  &  12.85     &   13.02     & 13.64    &  14.59 &   16.95  &  20.85\\                     
\multicolumn{1}{|c}{\multirow{1}{*}{}} & \multicolumn{1}{|c|}{\multirow{1}{*}{}}                                  &   $N=30$   &  13.01    &   13.17      & 13.84    &  14.77  &  17.11  &  20.72\\    
\cline{2-9}
\multicolumn{1}{|c}{\multirow{1}{*}{}} &   \multicolumn{2}{|c|}{Coverage Rate (\%)}                                & 100.00   & 100.00   & 100.00 & 99.98   & 98.99   & 85.56 \\
 \hline                                                                                                                                        
\end{tabular}
\end{center}
\caption{Comparison of the PTM and LWR models in estimating bin-based densities  using different sampling periods and penetration rates, for the period 5:15-5:30 pm. ``Coverage Rate" is the ratio of the number of active bins and the total number of bins.}
\label{tabPTMvsLWR3}
\end{table}

\subsection{Deviation of traffic from equilibrium}

In this section we calculate $q-q^*$, the indicator of the deviation of a vehicle speed from the equilibrium speed of surrounding vehicles. The PTM with a multi-valued FD allows disequilibrated traffic and complex waves, commonly observed  in congested highway traffic,  to be described and modeled. This is in contrast to the LWR model which captures only the equilibrium traffic states.

\begin{figure}[p!]
\includegraphics[width=1\textwidth]{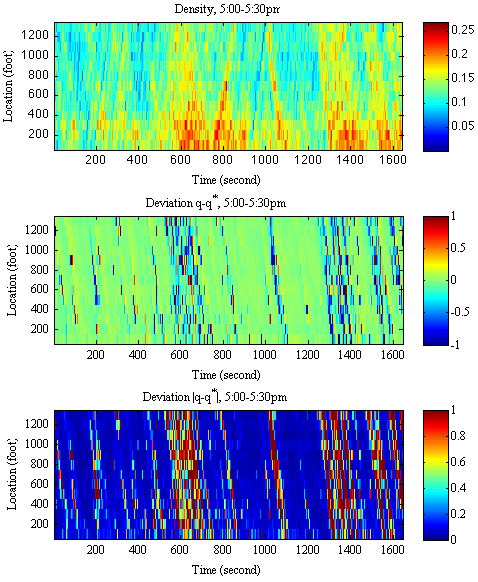}
\caption{Space-time diagram of the ground-truth vehicle density (top), signed deviation $q-q^*$ (center), and unsigned deviation $|q-q^*|$ (bottom). The deviations are calculated according to \eqref{eqn4}.}
\label{figqcombined}
\end{figure}

Using a similar procedure described in Section \ref{seccellrho}, we compute the quantity $\hat q=q-q^*$ using formula \eqref{eqn4}. It has been shown in Table \ref{tabssrhoqentireperiod} that the estimation of $\hat q$ is very susceptible to under sampling and, potentially, lower probe penetration as well, we utilize the entire ground-truth data ($N=1$, penetration rate $=100\%$) to reconstruct the bin-based $\hat q$, following a similar procedure described in Section \ref{seccellrho}.  The result is  shown in Figure \ref{figqcombined} for the traffic during 5:00-5:30pm. It can be seen from the top figure showing the density that the traffic was significantly congested and multiple shock waves were present and propagated with positive or negative speeds. The picture in the middle shows the bin-based $q-q^*$, which takes value between $[-1,\,1]$. We observe that traffic stayed close to equilibrium ($q-q^*=0$) where density was relatively low (typically below 0.16 vehicle/foot). However, significant deviation occurred with higher traffic densities. Furthermore, the presence of these high densities was accompanied by a mixture of above-equilibrium ($q-q^*>0$) and below-equilibrium ($q-q^*<0$) states.

We note that inside the same bin, above-equilibrium and below-equilibrium traffic may be present simultaneously; therefore, averaging $q-q^*$ within the same bin does not fully capture the deviation. Thus, we further compute the bin-based quantity $|q-q^*|\in[0,\,1]$, and the result is shown in the bottom picture of Figure \ref{figqcombined}. Similar to the signed deviation $q-q^*$, the unsigned deviation $|q-q^*|$ suggests a strong correlation with the vehicle density. This is consistent with the established observation that congested highway traffic tends to be unstable, due to inhomogeneous driving behavior \citep{LL} and the propagation of disturbances \citep{BABW}.

\section{Estimating instantaneous emission rate}\label{secHC}

It has been shown in Section \ref{secNumericalp} that the HC emission rate along vehicle trajectories are largely deteriorated by under sampling.  Such significant error is caused by the high variation in the acceleration profile. However, it is expected that when aggregated on a cell or link level, the estimation accuracy of these quantities may improve even with large sampling periods.  To confirm such an observation,  we will next perform Eulerian-based HC emission estimation.

\subsection{Bin-based estimation with insufficient data coverage}\label{subsecbinemission}

The bin-based emission rate can be estimated as follows. We first calculate the density $\rho_{ij}$ for bin $C_{ij}$ following the procedure illustrated in Section \ref{seccellrho} using Method 1.  Secondly, we use available vehicle trajectories to compute $r^{HC}(\tau_k,\,x(\tau_k))$ using formula \eqref{hcrate}.  Then all $r^{HC}(\tau_k,\,x(\tau_k))$ that fall within $C_{ij}$ are collected and averaged, which give rise to the average instantaneous HC emission rate per vehicle in this bin. Finally, this rate is multiplied by the number of vehicles present in this bin, which is given by the density, to yield the HC emission rate of $C_{ij}$.

Following this  procedure, we perform the Eulerian estimation of $r^{HC}$ for a combination of  lower sampling frequencies and penetration rates, which are then compared with the ground truth. The bin-based relative errors are averaged and shown in Table \ref{tabemissioncell} for the period 5:00-5:30pm. Compared to the Lagrangian (single-vehicle) estimation, the Eulerian estimation shows improved  accuracy due to the effect of space-time aggregation; the reader is referred to Table \ref{tabpowerdemandundersampling} for a comparison of results. Interestingly, the influence of under sampling reduces as the penetration rate becomes lower. This may be explained from the following two potential source of error:
\begin{enumerate}
\item  the lower sampling frequency that causes high variations in the emission rate to be ignored; 

\item the computational procedure which directly combines a microscopic emission model with a macroscopic traffic model, and creates certain mismatch in these very different modeling scales. 
\end{enumerate}
It is likely that with very high penetration rate (e.g. 100\%, 50\%), factor (1) is the dominating factor, while for lower penetration rates factor (2) (i.e. representing the emission rate inside a bin using a few microscopic vehicle trajectories) becomes the major source of error.

\begin{table}[h!]
\begin{center}
\begin{tabular}{cc|c|c|c|c|c|}  \cline{3-6}
  & & \multicolumn{4}{|c|}{Probe Vehicle Penetration Rate}
\\
\cline{2-6}
   & \multicolumn{1}{|c|}{\multirow{1}{*}{Sampling}} & \multicolumn{1}{|c|}{\multirow{2}{*}{100 \%}} & \multicolumn{1}{|c|}{\multirow{2}{*}{50 \%}} & \multicolumn{1}{|c|}{\multirow{2}{*}{20 \%}}  &  \multicolumn{1}{|c|}{\multirow{2}{*}{10 \%}} 
\\
 & \multicolumn{1}{|c|}{\multirow{1}{*}{Period}}  & \multicolumn{1}{|c|}{\multirow{2}{*}{}} & \multicolumn{1}{|c|}{\multirow{2}{*}{}} & \multicolumn{1}{|c|}{\multirow{2}{*}{}}  &  \multicolumn{1}{|c|}{\multirow{2}{*}{}} 
\\

\hline
  \multicolumn{1}{|c|}{\multirow{2}{*}{$r^{\text{HC}}$ Average}} & $N=10$        &  11.04         &14.09     &  19.36      &  23.67       \\
 \multicolumn{1}{|c|}{\multirow{2}{*}{Error (\%)}}                           &   $N=20$      &   16.66         &  17.96     &   20.52    &  23.21      \\   
 \multicolumn{1}{|c|}{\multirow{2}{*}}                                               &  $N=30$      &   19.55         &  20.40    &   21.75     &   23.20      \\       \cline{1-6}
   \multicolumn{2}{|c|}{Coverage Rate}                                                & 100.00 \%  & 100.00 \%  & 100.00 \%  & 99.88 \%    \\
 \hline

\end{tabular}
\end{center}
\caption{Estimation of bin-based hydrocarbon emission rate $r^{HC}$ using different sampling periods and penetration rates. The time period is 5:00-5:30pm. The estimation of vehicle density employs Method 1.}
\label{tabemissioncell}
\end{table}

\subsection{Correction factor approach for estimating aggregate emission rate}
\begin{figure}[h!]
\begin{minipage}[b]{.95\textwidth}
\centering
\includegraphics[width=1\textwidth]{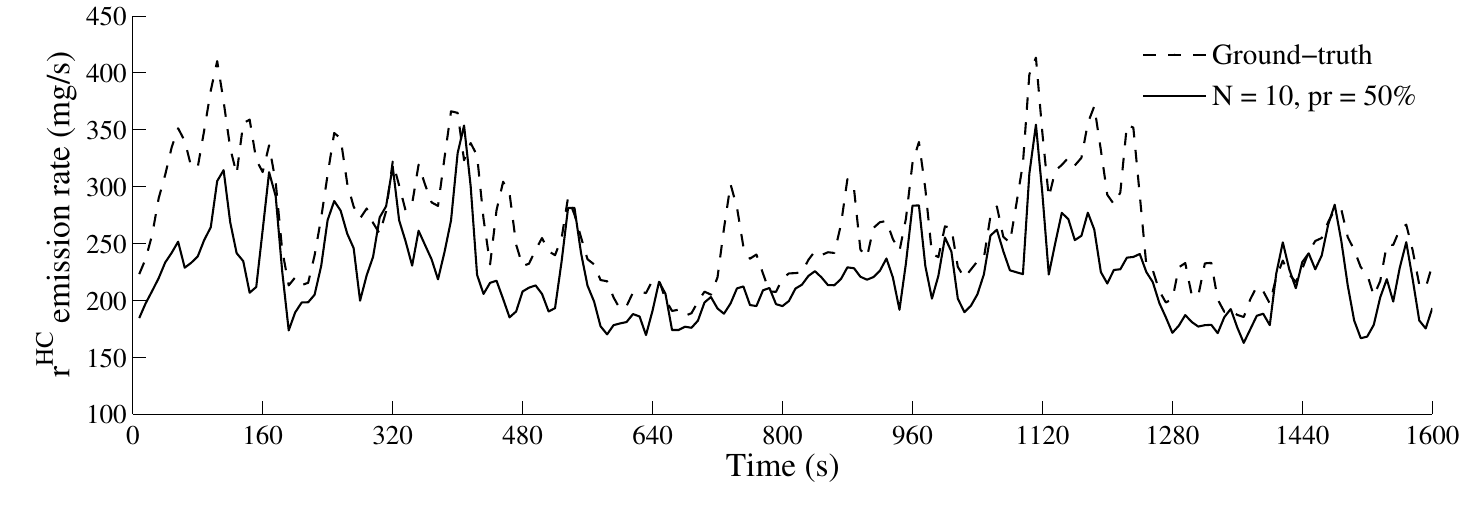}
\end{minipage}
\begin{minipage}[b]{.95\textwidth}
\centering
\includegraphics[width=1\textwidth]{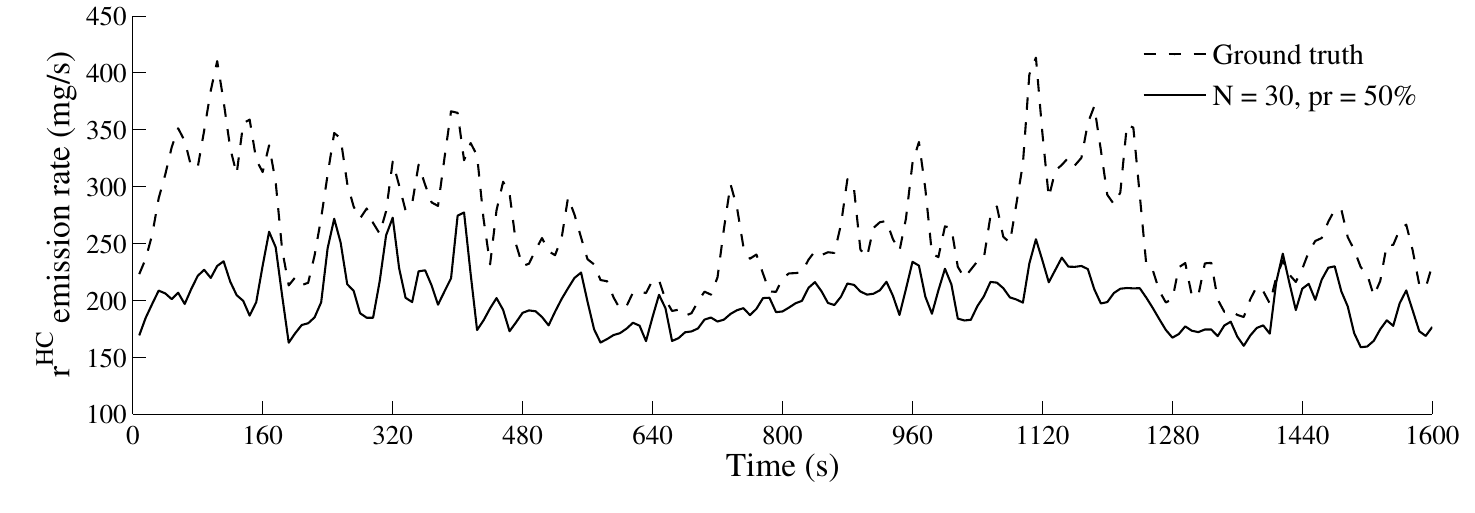}
\end{minipage}
\caption{Time-dependent HC emission rate on the entire study area (1600 feet in length) for the time period 5:00 - 5:30 pm. The top figure employs $N=10$ and 50\% penetration rate; the bottom figure employs $N=30$ and 50\% penetration rate.}
\label{figTogether_exp}
\end{figure}

To further aggregate the hydrocarbon emission rate, we compute the total emission rate along the entire road segment. This is done by summing up bin-based rates over all spatial cells for each time interval. An example of these time-dependent rates are shown in Figure \ref{figTogether_exp}, where the ground-truth $r^{HC}$ is compared with the estimated ones with [$N=10$, $pr=50\%$] and [$N=30$ and $pr=50\%$] respectively.  From these figures we see that the estimation captures to some degree the overall  trends of the emission rate, although it tends to underestimate the true value. The reason is that the under sampling overlooks higher-order variations in velocity and acceleration, and thus does not capture surges in the emission rate; a good example can be seen in Figure \ref{figEP}. And, such an underestimation is expected to become worse for larger sampling periods, as can be seen from a comparison between the two pictures in Figure \ref{figTogether_exp}.

The numerical result shown in Figure \ref{figTogether_exp} reveals discernible errors in emission estimation, although the overall time-varying trends of $r^{HC}$ is captured using our method. As we argued in Section \ref{subsecbinemission}, these errors are attributed to two factors, and to resolve these two issues simultaneously, we propose a correction factor approach to adjust the estimation process, which is demonstrably effective.

\subsubsection{Correction factors based on regression models}

\begin{figure}[h!]
\begin{minipage}[b]{.49\textwidth}
\centering
\includegraphics[width=1\textwidth]{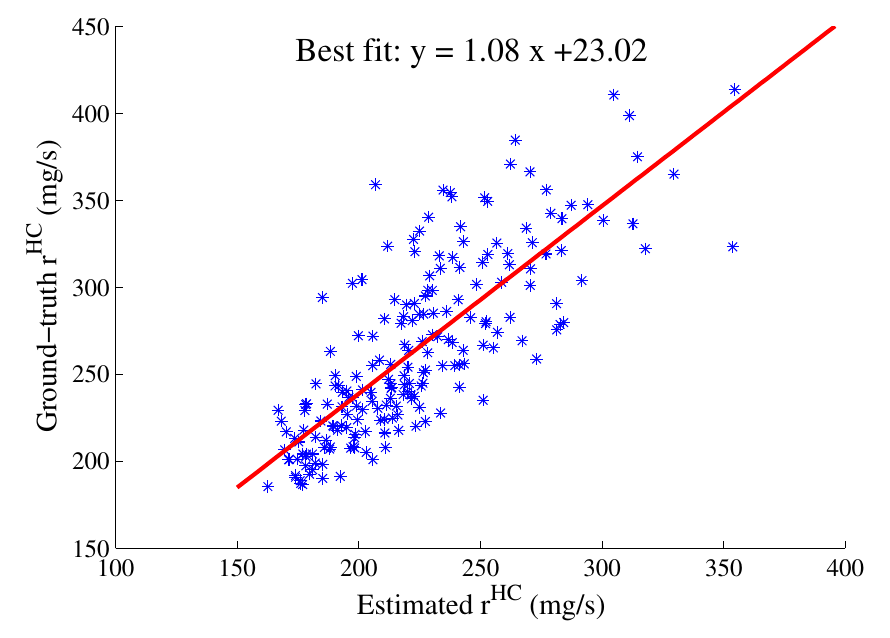}
\end{minipage}
\begin{minipage}[b]{.49\textwidth}
\centering
\includegraphics[width=1\textwidth]{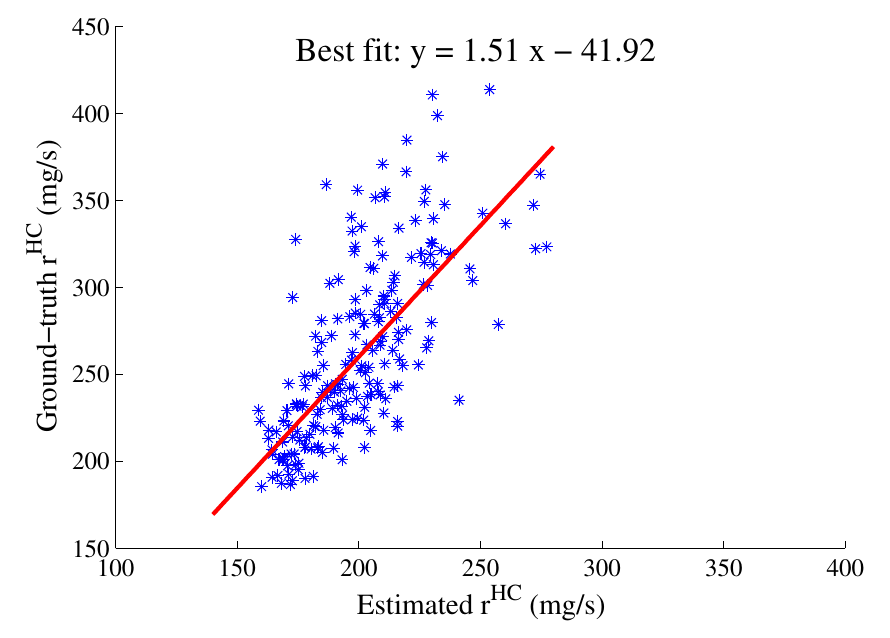}
\end{minipage}
\caption{Scatter plots of the estimated $r^{HC}$ vs. ground-truth $r^{HC}$ for the entire study area (1600 feet in length) during period 5:00 - 5:30 pm. The best linear fits are shown as solid lines. Left: $N=10$, $pr=50\%$; right: $N=30$, $pr=50\%$.}
\label{figlinearregression}
\end{figure}

We employ a linear regression approach which finds an appropriate affine relationship between the ground-truth data and the estimated ones.  To demonstrated the hypothesized affine relationship, we plot in Figure \ref{figlinearregression} the ground-truth value vs. the estimated value for the two cases illustrated in Figure \ref{figTogether_exp}. Each data point on the plane corresponds to one time interval;  the x-coordinate is the estimated value, and the y-coordinate is the true value. It is generally expected that different sampling periods and penetration rates require different linear regression coefficients. In fact, the best fit in the two cases  are: 
\begin{align}\label{correctionfactor1}
y~=~1.08 x +23.02 \qquad &(N=10, ~\hbox{pr}=50\%)
\\
\label{correctionfactor2}
y~=~1.51 x -41.92 \qquad &(N=30, ~\hbox{pr}=50\%)
\end{align}
Indeed, the $N=30$ case underestimate emission rates more than the $N=10$ case, and thus has a larger $1^{st}$-order coefficient.

\begin{figure}[h!]
\begin{minipage}[b]{.95\textwidth}
\centering
\includegraphics[width=1\textwidth]{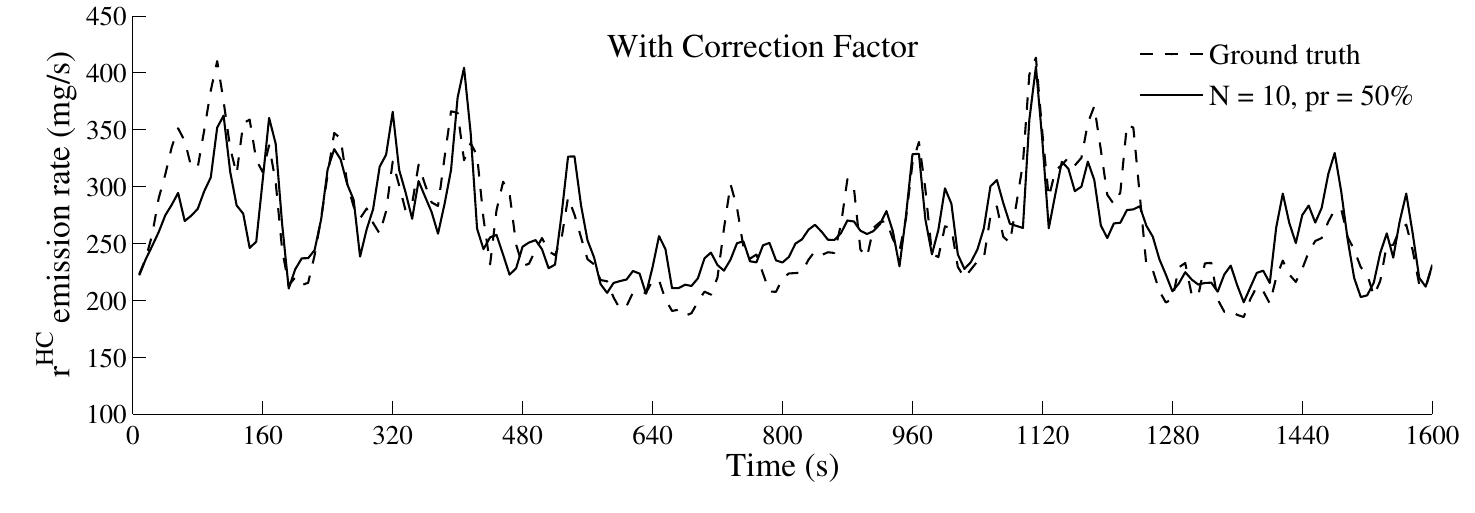}
\end{minipage}
\begin{minipage}[b]{.95\textwidth}
\centering
\includegraphics[width=1\textwidth]{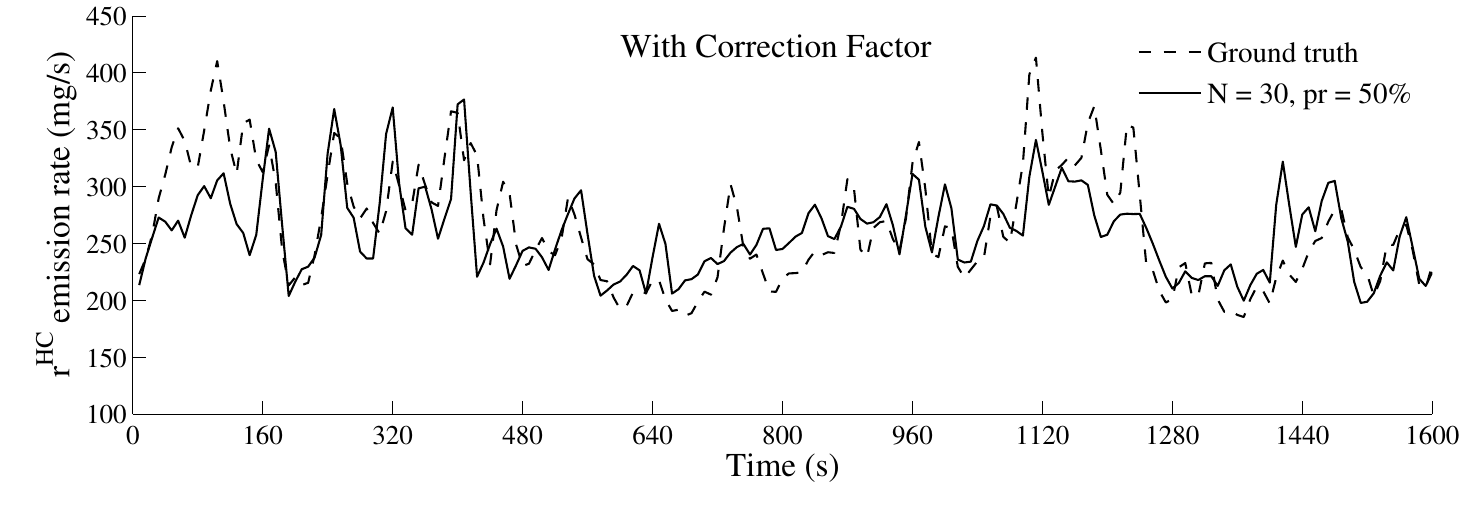}
\end{minipage}
\caption{The adjusted time series of $r^{HC}$ with the correction factors \eqref{correctionfactor1}-\eqref{correctionfactor2}.}
\label{figTogetherfitted}
\end{figure}

We denote by $r^{HC}_{true}$ the ground truth emission rate, and by $r^{HC}_{est}$ the estimated one. The linear model assumes that $r^{HC}_{true} = \beta_0+\beta_1 r^{HC}_{est}$ where $\beta_0,\,\beta_1\in \mathbb{R}$. Applying the affine adjustment \eqref{correctionfactor1} and \eqref{correctionfactor2}  yields the improved estimation depicted in Figure \ref{figTogetherfitted}. To rigorously justify the validity and robustness of the proposed correction factor approach, we employ a $k$-fold cross validation discussed below.

\subsubsection{Cross validation}

We partition available NGSIM vehicle trajectories into $k$ equal-size subsets. Among the $k$ subsets, a single subset is used as the training data to obtain the linear coefficients $\beta_0$ and $\beta_1$. The rest of the $k-1$ subsets are used to validate this affine adjustment by calculating the discrepancy between the predicted value and the ground-truth value. Such a procedure is repeated $k$ times, each with a distinct training dataset, and the errors from all $k$ repetitions are averaged. Such a process is known as $k$-fold cross validation.

We use the absolute relative error $|r_{est}^{HC}(t)-r_{true}^{HC}(t)| / |r_{true}^{HC}(t)|$ to indicate the accuracy of the estimation. A range of sampling periods and penetration rates are considered as shown in Table \ref{tabcv}, which summarizes results of the $k$-fold cross validation. Notice that for each penetration rate $pr$, the value $k$ is naturally chosen as $1/pr$.   It is quite clear from these cross validation results that applying the correction factor yields significantly improved estimation accuracy.

\begin{table}[h!]
\begin{center}
\begin{tabular}{cc|c|c|c|}  \cline{3-5}
  & & \multicolumn{3}{|c|}{Probe Vehicle Penetration Rate}
\\
\cline{2-5}
   & \multicolumn{1}{|c|}{\multirow{1}{*}{Sampling}} & \multicolumn{1}{|c|}{\multirow{2}{*}{50 \%}} & \multicolumn{1}{|c|}{\multirow{2}{*}{20 \%}}  &  \multicolumn{1}{|c|}{\multirow{2}{*}{10 \%}} 
\\
 & \multicolumn{1}{|c|}{\multirow{1}{*}{Period}}        & \multicolumn{1}{|c|}{\multirow{2}{*}{}}           & \multicolumn{1}{|c|}{\multirow{2}{*}{}}            &  \multicolumn{1}{|c|}{\multirow{2}{*}{}} 
\\

\hline
  \multicolumn{1}{|c|}{\multirow{2}{*}{$r^{\text{HC}}$ Average}} & $N=10$               & 10.34 (11.24)     &  12.34 (15.60)     &  13.47  (19.77)      \\
 \multicolumn{1}{|c|}{\multirow{2}{*}{Error (\%)}}                           &   $N=20$             &  12.21 (17.57)     &  13.65 (19.77)   &   14.32  (22.17)   \\   
 \multicolumn{1}{|c|}{\multirow{2}{*}}                                               &  $N=30$              &  13.52 (21.89)   &    14.31 (23.17)    &   14.39  (24.22)    \\       \cline{1-5}
  
 \hline                                
                     
\end{tabular}
\end{center}
\caption{Cross validation results of the correction factor approach. The errors in the parenthesis are from calculations without any correction factor.}
\label{tabcv}
\end{table}

\begin{remark}
It should be noted that the linear coefficients $\beta_0$ and $\beta_1$ employed by the correction factor approach are determined separately for each combination of sampling period and penetration rate. The reason is that these factors are devised to counter the effects of (1) underestimating the emission rates due to lower sampling frequency; and (2) possible mismatch between the microscopic emission model and macroscopic traffic model. Therefore, the factors must be determined separately as the effects (1) and (2) naturally vary in response to $N$ and the penetration rate.
\end{remark}

It can be seen that the estimation with correction factors has errors that are not significantly affected by under sampling or lower penetration rate, with a maximum difference of $4\%$ ($[N=10 ,\, pr=50\%]$ vs. $[N=30,\,pr=10\%$]) over all scenarios in Table \ref{tabcv}. This is in contrast to the case without the correction factors, which has a degradation of accuracy by as large as $13\%$ ($[N=10 ,\, pr=50\%]$ vs. $[N=30,\,pr=10\%$]).  This is due to the fact that the correction factors have taken into account misinterpretations of density and emission rate  caused by both larger $N$ and smaller $pr$, and hence the estimation is more resilient to the change in these parameters.

\section{Conclusion}\label{secconclusion}
This paper presents several approaches to estimate first- and higher-order traffic quantities based on  vehicle trajectory data. It also conducts extensive study of the accuracy of these estimations with respect to the probe penetration rate and sampling frequency of the mobile data. Most existing studies on mobile sensing tend to focus on the penetration rate, while the effect of  sampling frequency is rarely studied. We present several data fusion schemes that incorporate vehicle trajectory data into the second-order phase transition model, and perform the estimation of first-order quantities such as speed and density, and higher-order quantities such as acceleration, deviation (perturbation) and emission rate. The main findings are summarized below.
\begin{itemize}
\item Various Eulerian and Lagrangian  traffic quantities can be estimated via vehicle trajectory data, when the phase transition model (PTM) is employed.

\item For Lagrangian estimation, the first-order quantities are accurately reconstructed by the proposed method, even with low sampling frequency. For higher-order Lagrangian  quantities, lower sampling frequency significantly reduces the accuracy.  

\item The PTM-based estimation method is extended to perform Eulerian estimation,   with varying penetration rate and sampling frequency. The estimation accuracy shows a trend similar to that of the Lagrangian estimation, and deteriorates with lower penetration rate.

\item The PTM is compared with the LWR model in estimating Eulerian density with varying sampling frequency and penetration rate. The results show that the former yields more accurate results due to the fact that it  captures second-order deviation of traffic, which is frequently observed in congested and unstable traffic and is not described by first-order models.

\item A correction factor approach is applied to the Eulerian estimation of hydrocarbon emission rate, which significantly improves the accuracy. The method resolves the underestimation of emission due to insufficient sampling frequency, and the potentially biased emission estimation caused by mismatching modeling scales. As a result, the adjusted estimation is less sensitive to the changes in the sampling frequency and penetration rate. 

\end{itemize}

Further research will focus on applying the methodology to urban arterial networks and compare with the LWR-based approach. As we mentioned earlier in the introduction, the benefit of applying the PTM to arterial traffic will be somehow limited. However, higher-order variations in acceleration, emission and fuel consumption still need to be quantified and captured in connection with a given traffic flow model, and the effects of sampling frequency and penetration rate on these estimations should be properly understood in an urban environment.

\section{Acknowledgement} 
The authors would like to thank two anonymous reviewers for their constructive remarks.

\end{document}